\documentclass[12pt,reqno]{amsart}

\DeclareMathOperator{\id}{id}

\usepackage[margin=1in]{geometry}
\usepackage[T1]{fontenc}
\usepackage{lmodern}
\usepackage{microtype}
\usepackage{amsmath,amssymb,mathtools}
\usepackage{mathrsfs}
\usepackage{amsfonts}
\usepackage{xcolor}
\usepackage[shortlabels]{enumitem}
\usepackage[colorlinks,linkcolor=darkred,citecolor=darkblue,urlcolor=darkblue]{hyperref}
\usepackage{euflag}

\usepackage[normalem]{ulem}

\DeclareMathAlphabet{\lcal}{U}{dutchcal}{m}{n}

\definecolor{darkred}{rgb}{0.5,0,0}
\definecolor{darkblue}{rgb}{0,0,0.5}

\newtheorem{theorem}{Theorem}[section]
\newtheorem{corollary}[theorem]{Corollary}
\newtheorem{lemma}[theorem]{Lemma}

\newtheorem{remark}[theorem]{Remark}
\newtheorem{example}[theorem]{Example}

\numberwithin{equation}{section}
\allowdisplaybreaks

\newcommand{\C}{\mathbb C}
\newcommand{\Z}{\mathbb Z}
\newcommand{\N}{\mathbb N}
\newcommand{\cA}{\mathcal A}
\newcommand{\cB}{\mathcal B}
\newcommand{\cE}{\mathcal E}
\newcommand{\cF}{\mathcal F}
\newcommand{\cH}{\mathcal H}
\newcommand{\cK}{\mathcal K}
\newcommand{\cL}{\mathcal{L}}
\newcommand{\cM}{\mathcal M}
\newcommand{\cN}{\mathcal N}
\newcommand{\cO}{\mathcal O}

\newcommand{\cW}{\mathcal W}
\newcommand{\algA}{\mathfrak A}

\newcommand{\algAs}{\algA}
\newcommand{\spanop}{\operatorname{span}}
\newcommand{\ext}{\operatorname{ext}}

\newcommand{\vd}{{\tt{d}}}

\newcommand{\otn}{\otimes_{\min}}
\newcommand{\otm}{\otimes_{\max}}

\newcommand{\idA}{\mathtt{1}}
\newcommand{\idY}{\text{e}}
\newcommand{\idGroupAlg}{\idA}

\newcommand{\groupY}{\mathfrak{Y}}
\newcommand{\groupS}{\mathfrak{H}}

\newcommand{\df}[1]{{\textit{#1}}{\index{#1}}}

\newcommand{\simu}{\sim_{\text{u}}}

\newcommand{\lly}{\mathcal{Y}} %

\newcommand{\HSH}{\mathcal{H}}
\newcommand{\HSE}{\mathcal{E}}

\newcommand{\lfrac}{\operatorname{\ell\text{-}frac}}

\newcommand{\wtau}{\widetilde{\tau}}

\newcommand{\kk}{k}

 \title[Fej\'er--Riesz factorization for free products]{Fej\'er--Riesz factorization for free products}

\author[A.\ Jindal]{Abhay Jindal${}^{1,Q}$}
\address{Abhay Jindal, Faculty of Mathematics and Physics, University of Ljubljana, Slovenia}
\email{abhay.jindal@fmf.uni-lj.si}
\thanks{${}^1$Supported by the Slovenian Research Agency 
program P1-0222 and grant J1-50002.}

\author[I.\ Klep]{Igor Klep${}^{2,Q,R}$}
\address{Igor Klep, Faculty of Mathematics and Physics, University of Ljubljana 
\& Famnit, University of Primorska, Koper 
\& Institute of Mathematics, Physics and Mechanics,
Ljubljana, Slovenia}
\email{igor.klep@fmf.uni-lj.si}
\thanks{${}^2$Supported by the Slovenian Research Agency 
	program P1-0222 and grants J1-50002,
	N1-0217,
J1-60011, J1-50001, J1-3004 and J1-60025.}

\author[J. Levenson]{Jacob Levenson}
\address{Jacob Levenson, Department of Mathematics\\
	University of Florida\\ Gainesville } %
\email{levenson.j@ufl.edu}

\author[S. McCullough]{Scott McCullough}
\address{Scott McCullough, Department of Mathematics\\
	University of Florida\\ Gainesville} %
    
\email{sam@math.ufl.edu}

\thanks{${}^Q$This work was performed within the project COMPUTE, funded within the QuantERA II 
Programme that has received funding from the EU's H2020 research and innovation programme under the GA No.~101017733 {\normalsize\euflag}}

\thanks{${}^R$Part of this work was performed within the project SDPcode, funded within the QuantERA III 
Programme that has received funding from the EU's Horizon research and innovation programme under the GA No.~101212998 {\normalsize\euflag}}



\date{}

\begin{document} 

\begin{abstract}
We establish an operator-valued Fej\'er--Riesz factorization
theorem for the algebraic free product
$\algA=\C[F_r]*\cA_1*\cdots *\cA_s$,
where $F_r$ is a free group and each $\cA_j$ is a
finite-dimensional unital $C^*$-algebra.
With respect to a natural reduced-word length, called
extent, we show that every element $A\in B(\cH)\otimes\algA$
of extent at most $2d$ that is positive under all
finite-dimensional unital $*$-representations of $\algA$
admits a factorization $A=B^*B$, with $B$ of extent
at most $d$.
This sharp half-extent bound
holds without strict positivity.
The proof combines a new finite-dimensional realization
theorem for positive truncated kernels with
operator-system methods.
We also establish factorization results for tensor
products with group algebras of arbitrary discrete groups
under strict positivity in the maximal $C^*$-tensor product;
for amenable groups, strict positivity in the minimal
tensor product suffices.
\end{abstract}

\maketitle

\tableofcontents

\section{Introduction}

The classical Fej\'er--Riesz theorem asserts that a nonnegative trigonometric polynomial on the unit circle admits a square factorization. In operator-valued form, if
\begin{equation}
    \label{e:classic-fr}
 p(e^{i\theta})
 \,=\,
 \sum_{n=-k}^k A_n e^{in\theta},
 \qquad A_n\in B(\cH),
\end{equation}
is positive semidefinite for every $\theta$, then there is an analytic polynomial
$q(e^{i\theta})=\sum_{n=0}^k B_n e^{in\theta}$ with coefficients in $B(\cH)$ such that $p=q^*q$; see \cite{Ro, DR10}. In real algebraic geometry \cite{BCR98,Mar08,Sc24}, the result is understood as a sums-of-squares representation or a Positivstellensatz. 
Beyond positivity certificates, Fej\'er--Riesz factorization is closely tied to the operator-algebraic theory of Toeplitz systems \cite{CvS21} and has concrete applications in 
prediction theory \cite{Rose},
digital filter design \cite{WBV99,Dum17}, and  quantum signal processing \cite{BS25}. 

Generalizations include multivariate Fej\'er--Riesz theorems
\cite{GW04,Dritschel,DW05} and factorization and extension results
for noncommutative polynomials and free groups
\cite{Mc,He,BT07,Oza13}. Recent group-algebra results concern
free products of copies of $\Z_2$ \cite{KLM} and free products
of finite abelian groups \cite{JKM}. The latter work introduced
the extent filtration as a degree-like measure for reduced words.

In this paper we consider the algebraic free product
\[
 \algA=\C[F_r]*\cA_1*\cdots *\cA_s,
\]
where $F_r$ is a free group and each $\cA_j$ is a
finite-dimensional unital $C^*$-algebra. Our main result
establishes an exact Fej\'er--Riesz factorization with
operator coefficients and the sharp half-extent bound,
without assuming strict positivity. 
In particular, it covers
free products with matrix-algebra factors and with group
algebras of non-abelian finite groups.

A central ingredient here is a new realization theorem for
positive truncated kernels, extending local left-multiplication
operators to representations of the free-product factors.
Together with operator-system methods, this yields the
factorization and the stated extent bound.\looseness=-1

We also study tensor products with group algebras of arbitrary
discrete groups. In that setting, strict positivity in the
maximal $C^*$-tensor product yields a factorization with the
same half-extent bound. The maximal positivity assumption
cannot in general be replaced by positivity in spatial
representations; for amenable groups, however, strict
positivity in the minimal tensor product suffices.

We next recall the free-product notation and the extent
filtration needed to state the results precisely.

\subsection{The algebraic free product}
\label{ss:algA}
Let \df{$F_r$} denote the free group on the $r$ generators $x_1,\dots,x_r$ and let \df{$\C[F_r]$} denote the resulting group algebra:
\[
\C[F_r] 
\,\cong\,
\C[x_1,x_1^{-1}]*\cdots *\C[x_r,x_r^{-1}].
\]
For $1\le j\le s$, let \df{$\cA_j$} be a
finite-dimensional unital $C^*$-algebra, and let
\begin{equation}\label{eq:free-product-algebra}
\algA
\,=\,
\C[F_r]*\cA_1*\cdots *\cA_s
\end{equation}
be the algebraic unital free product. \index{$\algA$}
For each infinite cyclic factor $\C[x_i,x_i^{-1}]$, let \df{$\tau_{x_i}$} denote the \df{canonical trace},
\[
\tau_{x_i}(x_i^n)=
\begin{cases}
	1, & n=0,\\
	0, & n\neq 0.
\end{cases}
\]
For each finite-dimensional factor $\cA_j$, fix a faithful tracial state \df{$\tau_j$} and write
\[
\cA_j=\C \idA \oplus \cA_j^\circ,
\qquad
\cA_j^\circ=\{a\in \cA_j:\tau_j(a)=0\}.
\]
Choose an orthonormal basis \df{$\cB_j$} of \df{$\cA_j^\circ$} with respect to the inner product
\[
\langle a,b\rangle_j=\tau_j(b^*a).
\]

\subsection{Reduced words}
\label{ss:reduce}
A \df{reduced word} in $\algA$ is either the identity \df{$\idA$} or a product
\[
w=w_1w_2\cdots w_m
\]
such that each syllable $w_\ell$ is either a nontrivial power $x_i^n$, with $n\in\Z\setminus\{0\}$, or an element of one of the bases $\cB_j$, and consecutive syllables come from different free product factors. Thus two consecutive powers of the same generator $x_i$ are not allowed, and two consecutive letters from the same finite-dimensional factor are not allowed. Let \df{$\cW$} denote the set of reduced words in $\algA.$

\smallskip

The \df{free-product trace} \index{$\tau$}
\[
\tau
\,=\,
\tau_{x_1}*\cdots *\tau_{x_r}*\tau_1*\cdots *\tau_s
\]
is the linear functional on $\algA$ that extracts the coefficient of the empty word. Equivalently, $\tau(\idA)=1$ and $\tau(w)=0$ for every nonempty reduced word $w$.
We refer the reader to \cite{VDN92} for a detailed discussion of the free product of $C^{*}$-algebras. 

\smallskip

\subsection{Extent}
\label{ss:extent}
Define the \df{weight} of a syllable by
\[
\operatorname{wt}(w_\ell)=
\begin{cases}
	|n|, & w_\ell=x_i^n,\ n\neq 0,\\[1mm]
	1, & w_\ell\in \cB_j.
\end{cases}
\]
The \df{extent} of a reduced word $w=w_1\cdots w_m$ is
\[
\ext(w) 
\,=\,
\sum_{\ell=1}^m \operatorname{wt}(w_\ell),
\qquad
\ext(\idA)=0.
\]
Let \df{$\cW_\vd$} be the finite set of reduced words of extent at most $\vd$ and set \index{$\algAs_\vd$} \index{$N(\vd)$}
\[
\algAs_\vd=\spanop \cW_\vd\subseteq \algA.
\]

\subsection{Involution}
\label{ss:involute}
The involution on $\algA$ is induced by
$(x_i^n)^*=x_i^{-n}$ and the involutions on the factors
$\cA_j$, with the order of syllables reversed.
The adjoint of a basis element of $\cA_j^\circ$ need not
belong to the chosen basis $\cB_j$, so the adjoint of a
reduced word may be a linear combination of reduced words.
Nevertheless, the spaces $\algA_\vd$ are invariant under
involution; see Lemma~\ref{lem:adjoint-invariant}.
For operator-valued polynomials we use
$(T\otimes a)^*=T^*\otimes a^*$ and extend conjugate-linearly.

\subsection{The main result}
\label{ss:thm-main}
Let $\cH$ be a Hilbert space and let $A\in B(\cH)\otimes \algAs_{2\vd}$. Writing
\[
A 
\,=\,
\sum_{\ell=1}^{N(2\vd)} A_\ell\otimes r_\ell,
\qquad A_\ell\in B(\cH),\quad r_\ell\in\cW_{2\vd},
\]
we say that $A$ is \df{positive  on finite-dimensional representations}, and write $A\succeq 0$, if $A=A^*$ and for every finite-dimensional Hilbert space $\cF$ and every unital $*$-representation
$ \sigma:\algA\to B(\cF),$ %
\[
(\id_{B(\cH)}\otimes \sigma)(A) 
\,:=\,
\sum\limits_{\ell =1}^{N(2\vd)} A_{\ell} \otimes \sigma (r_{\ell})
\succeq 0
\quad\text{on }\cH\otimes \cF.
\]

Our first  main result is the following.

\begin{theorem}\label{thm:main}
Let $\vd$ denote a positive integer and let $\cH$ be a Hilbert space.  If
$  A\in B(\cH)\otimes \algAs_{2 \vd} $
is positive on every finite-dimensional unital $*$-representation of $\algA$, then there exists a Hilbert space $\cE$ and an element
$ B 
\in B(\cH,\cE)\otimes \algAs_\vd $
such that
\[
   A=B^*B.
\]
We can take $\cE=\cH\otimes \C^{N(\vd)}$, where $N(\vd)= {\rm dim}(\algA_{\vd})$.
\end{theorem}

The converse of Theorem~\ref{thm:main} is immediate: a factorization $A=B^*B$
certifies positivity in every representation. The theorem
requires no strict positivity and gives the uniform
ampliation bound $\cE=\cH\otimes\C^{N(\vd)}$.

The free-group case ($s=0$) has antecedents in
\cite{Mc,BT07,Oza13}; see, e.g.,
\cite[Theorem~26]{Oza13} and \cite[Theorem~1.3]{JKM-ieot}
for operator-valued factorizations with prescribed support.
For free products of finite abelian groups,
\cite[Theorem~B]{JKM} gives factors of extent at most
$\vd+1$ for an input of extent at most $2\vd$.
Theorem~\ref{thm:main} improves this general bound to
$\vd$ and allows arbitrary finite-dimensional
$C^*$-algebra factors.

\smallskip

As a consequence of Theorem~\ref{thm:main}, we obtain the corresponding
Fejér--Riesz factorization on free products
of a free group and finitely many finite groups.

\begin{corollary}\label{cor:group-version}
Let
	\begin{equation}\label{eq:ourgroup}
	\mathbb G
	\,=\,
	F_{r}
	*
	\mathbb G_1*\cdots * \mathbb G_s,
	\end{equation}
where each $\mathbb G_j$ is a finite group.
Equip $\mathbb G$ with the notion of extent defined above,
using $\mathbb G_j\setminus\{\idA\}$ in place of
$\mathcal B_j$. If
    
    \[
	A \,=\, \sum_{ \ext(g) \leq  2\vd} A_g\otimes g
	\in B(\mathcal H)\otimes \mathbb C[\mathbb G] 
    \]
	is a self-adjoint operator-valued group polynomial of extent at most $2\vd$ such that for every finite-dimensional Hilbert space $\mathcal F$ and
	every unitary representation
	$\pi: \mathbb G\to B(\mathcal F)$
	\[
	(\operatorname{id}_{B(\mathcal H)}\otimes \pi)(A) \succeq  0,
	\]
	then there exists
  $B = \sum \limits_{ \normalfont \ext(g) \leq \vd}  B_g \otimes g \in B(\HSH, \mathcal H \otimes \C^{N(\vd)}) \otimes \mathbb{C}[\mathbb{G}] $  
	such that
	\[
	A=B^*B.
	\] 
\end{corollary}

The non-quantitative scalar sums-of-squares conclusion for the groups in Corollary \ref{cor:group-version} is already known. Indeed, a group of the form \eqref{eq:ourgroup}  is virtually free. Netzer–Thom \cite[Theorem 6.1]{NT13} establish the scalar free-group result and remark its extension to virtually free groups, with the proof of the latter later given in \cite[Theorem 3.2.1(b)]{GC23}. These results cover the qualitative scalar conclusion, but do not yield either the prescribed extent bound or the operator-valued factorization.
Thus Corollary \ref{cor:group-version} strengthens these conclusions by allowing operator coefficients, and by guaranteeing that an element of extent at most $2\vd$ admits a factorization with factors of extent at most $\vd$. 

\begin{remark}\rm
 \label{r:rep-bound}
  In Theorem~\ref{thm:main} and Corollary~\ref{cor:group-version}, in  the case that $\cH$ has finite dimension $K,$ it suffices to check that $A$ is positive against unital $\ast$-representations
of $\algA$ on a Hilbert space of dimension at most
\[
\operatorname{lcm}
\bigl\{
1,2,\ldots, K\, N(\vd), 
\dim\cA_1,\ldots,\dim\cA_s
\bigr\}.
\]
\end{remark}

\subsection{Tensoring with a group algebra}
\label{ss:add-Y}

Let $\groupS$ be a discrete group. Given a 
subset $\mathbf y\subseteq\groupS$, let $\lly$
denote the monoid generated by $\mathbf y$.
For $M\ge0$, let $\lly_{\le M}$ consist of
the elements of $\lly$ representable by words
in $\mathbf y$ of length at most $M$.
For a subset $E$ of a group $G$, let
\[
 \lfrac E=\{u^{-1}v:u,v\in E\}.
\]

Our next theorem gives an operator-valued factorization
for arbitrary additional groups under maximal
$C^*$-tensor-product positivity. Here $C_u^*(\algA)$
denotes the universal $C^*$-completion of $\algA$, and
$C^*(\groupS)$ is the full group $C^*$-algebra of $\groupS$.

\begin{theorem}
\label{t:main+maxY}
Let $\groupS$ be an arbitrary discrete group, let $\vd$ be a
positive integer, and let $\cH$ be a Hilbert space.
Suppose that
\[
 A=A^*\in B(\cH)\otimes\algA_{2\vd}\otimes\C[\groupS]
\]
and
\begin{equation}\label{e:maxY-positive}
 A\succ0
 \quad\text{in}\quad
 B(\cH)\otm C_u^*(\algA)\otm C^*(\groupS).
\end{equation}
(Here $A\succ0$ means that $A$ is positive and
invertible in the indicated $C^*$-algebra.)

Then there exist a positive integer $n$, a finite subset
$F\subseteq\groupS$, and operators
\[
 B_{q,g}\in B(\cH,\cH^{\oplus n}),
 \qquad q\in\cW_\vd,\quad g\in F,
\]
such that
\[
 B=\sum_{\substack{q\in\cW_\vd\\g\in F}}
       B_{q,g}\otimes(q\otimes g)
 \qquad\text{satisfies}\qquad A=B^*B.
\]
In particular, the target space $\cE=\cH^{\oplus n}$ is
separable whenever $\cH$ is separable.

Suppose in addition that $\groupS=\lfrac\lly$, where
$\lly$ is the monoid generated by a finite set
$\mathbf y\subseteq\groupS$, and that $A$ has the form
\[
 A=\sum_{\substack{a\in\cW_{2\vd}\\
                      g\in\lfrac\lly_{\le M}}}
 A_{a,g}\otimes(a\otimes g),
 \qquad A_{a,g}\in B(\cH),
\]
for some positive integer $M$. Then $F$ can be chosen to
be contained in $\lly_{\le M'}$ for some $M'\ge M$.
\end{theorem}

Since the group support of $A$ is finite,
Lemma~\ref{r:reduce-to-finite} shows that, in the
first assertion of Theorem~\ref{t:main+maxY},
we may assume without loss of generality that
$\groupS$ is finitely generated.

For amenable supporting subgroups, the maximal positivity hypothesis
can be expressed in terms of spatial evaluations. This gives the
following corollary of Theorem~\ref{t:main+maxY}.

\begin{corollary}
\label{t:main+amenableY}
Let $\groupS$ be an amenable group, and let
$\lly$ be the monoid generated by a finite set
$\mathbf y\subseteq\groupS$. Suppose that
\[
 \groupS=\lfrac\lly.
\]
Let $\vd,M$ be positive integers and let $\cH$ be a
Hilbert space. Suppose that $A=A^*$ has the form
\[
 A=
 \sum_{\substack{a\in\cW_{2\vd}\\
                 g\in\lfrac\lly_{\le M}}}
 A_{a,g}\otimes(a\otimes g),
 \qquad A_{a,g}\in B(\cH).
\]
If 
\begin{equation}\label{eq:positive-evaluations}
	A\succ0
 \quad\text{in}\quad
 B(\cH)\otn C_u^*(\algA)\otn C^*(\groupS),
\end{equation}
 then there
exist an integer $M'\ge M$, a Hilbert space $\cE$,
and operators
\[
 B_{q,u}\in B(\cH,\cE),
 \qquad q\in\cW_\vd,\quad u\in\lly_{\le M'},
\]
such that
\[
 B=
 \sum_{\substack{q\in\cW_\vd\\u\in\lly_{\le M'}}}
 B_{q,u}\otimes(q\otimes u)
 \qquad\text{satisfies}\qquad A=B^*B.
\]
The space $\cE$ can be chosen to be a finite ampliation of
$\cH$. Thus $\cE$ can be chosen separable whenever $\cH$
is separable. If $\cH$ is infinite-dimensional, the
coefficients of $B$ may be taken in $B(\cH)$.\looseness=-1
\end{corollary}

{The strict positivity condition in equation~\eqref{eq:positive-evaluations} is equivalent to the existence
of an  $\varepsilon>0$ such that
\begin{equation*} %
 (\id_{B(\cH)}\otimes\pi)(A)\succeq\varepsilon I
\end{equation*}
for every unital $*$-representation $\pi$ of
$\algA\otimes\C[\groupS]$ on a Hilbert space.}

For matrix coefficients, the amenability assumption can be omitted.

\begin{remark}\rm
\label{r:matrix-coefficients}
For $\cH=\C^K$, the conclusion of
Corollary~\ref{t:main+amenableY} holds without assuming
$\groupS$ amenable if \eqref{eq:positive-evaluations} is replaced by
$(\id_{M_K}\otimes\pi)(A)\succ0$ for every unital
$*$-representation $\pi$ of
$\algA\otimes\C[\groupS]$.
Indeed, put $D=C_u^*(\algA)\otm C^*(\groupS)$.
If $A$ were not strictly positive in $M_K(D)$, a state
attaining the minimum of its spectrum would yield a GNS
representation of $M_K(D)$ in which $A$ is not strictly
positive. Since every unital representation of $M_K(D)$
is of the form $\id_{M_K}\otimes\pi$, up to unitary
equivalence, this contradicts the hypothesis.
Thus $A\succ0$ in $M_K(D)=M_K(\C)\otm D$,
and Theorem~\ref{t:main+maxY} applies.
As $\groupS$ is finitely generated, separable
representations suffice.
\end{remark}

Strict positivity in Remark~\ref{r:matrix-coefficients} is essential in general,
as the conclusion can fail even for scalar-valued nonnegative
trigonometric polynomials in three commuting variables
(cf.~\cite{Sc00}). 
Example~\ref{ex:below} shows that, even for an amenable
additional group, strict positivity of the quadratic form
in every spatial evaluation does not suffice for
$B(\cH)$-valued coefficients.
Even a uniform lower bound in all spatial evaluations need not suffice
for general additional groups, as Example~\ref{r:operator-coefficients-obstruction}
in Section~\ref{s:operator-coefficients} shows. This explains the role
of maximal tensor-product positivity in Theorem~\ref{t:main+maxY}.
Corollary~\ref{c:main+finiteY} below illustrates how the positivity
hypothesis can be relaxed when the additional group $\groupS$ is finite.

Remark \ref{r:matrix-coefficients} is related to \cite[Theorem 1.1]{KLM}, which already establishes a strict matrix-valued mixed-group factorization with preservation of support in a free coordinate. \cite{KLM} also treats non-strict operator-valued factorization for free products of copies of $\mathbb Z_2$ \cite[Theorem~1.6]{KLM}; see also \cite[Remark~1.7]{KLM}. The support conventions, however, are different: \cite{KLM} employs shortlex and analytic/hereditary support restrictions, whereas the present results use extent and allow inverses in the free-group coordinate. (Compare Examples 8.1 and 8.2 of \cite{KLM}). 

Unlike the compactness argument of \cite[Sections~6--7]{KLM}, the matrix-coefficient assertion of Remark~\ref{r:matrix-coefficients} follows directly from Theorem~\ref{t:main+maxY}, including the extent bound $\vd$ for the present class of free products.

When $\groupS$ is finite, Theorem \ref{thm:main} yields the following stronger conclusion.

 \begin{corollary}
\label{c:main+finiteY}
   Let $\groupS$ be a finite group, let $\vd$ be a positive
integer, and let $\cH$ be a Hilbert space.
If
\[
 A\in B(\cH)\otimes\algA_{2\vd}\otimes\C[\groupS]
\]
is nonnegative against every 
unital $*$-representation
of $\algA\otimes\C[\groupS]$ on a finite-dimensional Hilbert space,
then there exists
\[
 B\in B(\cH,\cE)\otimes\algA_\vd\otimes\C[\groupS],
 \qquad
 \cE=\cH\otimes\C^{N(\vd)},
\]
such that $A=B^*B$.
 
 Further, if $\cH$ is finite-dimensional, then it suffices to verify that $A$ is positive %
against unital $*$-representations of
$\algA\otimes\C[\groupS]$ on Hilbert spaces of dimension
at most
 \[
\operatorname{lcm}
\bigl\{
1,2,\ldots,\dim \cH\, |\groupS|\, N(\vd),
\dim\cA_1,\ldots,\dim\cA_s
\bigr\}.
\]
\end{corollary}

\subsection{Comparison with earlier approaches}
\label{ss:earlier-approaches}

The approach of \cite{JKM} to free products of finite
abelian groups uses a Positivstellensatz for positive
operator-valued measures and Boca's theorem on free
products of completely positive maps.
It exploits the commutative spectral structure of
finite abelian group algebras, which is absent for
general finite-dimensional $C^*$-algebras.
The methods of \cite{KLM} rely on positive semidefinite
matrix completions adapted to the shortlex order
on free semigroups and free products of copies of
$\Z_2$. These constructions do not directly yield
representations respecting the multiplication
relations within arbitrary finite-dimensional
$C^*$-algebra factors.

{The principal novelty here underlying Theorems~\ref{thm:main} and \ref{t:main+maxY} and their corollaries is the realization
of positive truncated kernels by representations
respecting the relations of all free-product factors,
while retaining control of the extent.  See Theorem~\ref{prop:main:rep},  which extends
local left-multiplication operators to representations
of the free-product factors.
Combined with operator-system methods, this yields Theorem~\ref{thm:main}, including
the sharp half-extent factorization without strict
positivity.
}

\subsection{Reader's guide}\label{ss:guide}
  In outline, the proof of Theorem~\ref{thm:main} proceeds as follows.  
The positivity of $A$ is interpreted as complete
positivity of a linear map $\Phi_A$ defined on an
operator system $\cO_\vd$. When
$\algA=\C[x,x^{-1}]$, this operator system is,
after ordering the basis appropriately, the space of
$(2\vd+1)\times(2\vd+1)$ Toeplitz matrices.
Positive elements of $\cO_{\vd}  \otimes M_K(\C),$ where $M_K(\C)$ is the space of complex $K\times K$ matrices,  
are seen, via a finite-dimensional dilation result for positive truncated kernels (Theorem \ref{prop:main:rep}),
   to arise from representations of $\algA.$  Applying the Arveson Extension Theorem to $\Phi_A$ and factoring the resulting positive semidefinite Choi matrix associated to the extended version of $\Phi_A$ completes the proof.

For Theorem~\ref{t:main+maxY}, we use an algebraic
separation argument, in the spirit of \cite{HM04},
and a group-equivariant form of the truncated-kernel
realization. If $A - \delta I$ is not a sum of squares of elements of $\algA$-extent at most $\vd$, an algebraic order-unit separation argument produces a separating linear functional $\Lambda$.     
Applying Theorem~\ref{prop:main:rep} with
$G=\mathcal U(B(\cH))\times\groupS$
realizes $\Lambda$ as a vector functional on a
representation of the maximal $C^*$-tensor
product, contradicting the strict positivity of $A$.
Combining Theorem~\ref{t:main+maxY} with Pisier's free-product theorem, Lemma~\ref{lem:min-max-arbitrary-H}, and the nuclearity of $C^*(\groupS)$ for amenable $\groupS$ yields Corollary~\ref{t:main+amenableY}.    %

Section~\ref{s:prelims} collects some needed facts about
$\algA$. Section~\ref{s:kernels-and-reps} develops the
representation theory of positive truncated kernels
and proves 
Theorem \ref{prop:main:rep}.
Section~\ref{s:proofs} establishes Theorem~\ref{thm:main}
and the associated free-product and finite-group Corollaries~\ref{cor:group-version} and \ref{c:main+finiteY}.
Section~\ref{s:operator-coefficients} proves the results on
tensor products with group algebras (Theorem~\ref{t:main+maxY} and Corollary~\ref{t:main+amenableY}), and gives an obstruction
for general groups in Example~\ref{r:operator-coefficients-obstruction}.

\subsection{Use of AI}  
The authors acknowledge the use of LLMs during the late stages of work on this paper.  It responded to our prompts to unify two instances of the matrix-coefficient factorization now stated in Remark~\ref{r:matrix-coefficients} - one without the free group factor in $\algA$ and general $\groupS$ and the other with general $\algA$ and $\groupS$ an amenable group - using Lemma~\ref{l:key-extend}. It also produced Example~\ref{r:operator-coefficients-obstruction} and was used to assist in final editing. All mathematical details were independently checked and verified by the authors, who assume full responsibility for their mathematical correctness.

\section{Preliminaries}\label{s:prelims}

\begin{lemma}\label{lem:orthonormal-reduced-words}
The reduced words in $\algA$ form an orthonormal basis of the pre-Hilbert  space obtained from $\algA$ with the inner product
\[
   \langle u,v\rangle=\tau(v^*u).
\]
\end{lemma}

\begin{proof}
For each cyclic factor put
\[
   \C[x_i,x_i^{-1}]^\circ=\spanop\{x_i^n:n\neq 0\}.
\]
As a vector space, the algebraic free product decomposes as
\[
\algA
\,=\,
\C \idA\oplus
 \bigoplus_{m\geq 1}
 \bigoplus_{\ell_j\ne \ell_{j+1} \, (1\le j <m)}  %
 E_{\ell_1}^\circ\otimes E_{\ell_2}^\circ\otimes\cdots\otimes E_{\ell_m}^\circ,
\]
where  either $E_\ell^\circ =\C[x_i,x_i^{-1}]^\circ $ or $E_\ell^\circ =\mathcal{A}_j^\circ$. 
The free-product inner product is the tensor-product inner product on each summand, and distinct summands are  orthogonal to one another. Since
\[
   \{x_i^n:n\neq 0\}
\]
is an orthonormal basis for $\C[x_i,x_i^{-1}]^\circ$ and $\cB_j$ is an orthonormal basis for $\cA_j^\circ$, the corresponding reduced words form an orthonormal basis. 
\end{proof}

\begin{lemma}\label{lem:adjoint-invariant}
For every nonnegative integer $\vd$, the subspace $\algAs_\vd$ is closed under the involution.
\end{lemma}

\begin{proof}
It suffices to consider a reduced word $w=w_1\cdots w_m\in\algAs_\vd$. If $w_\ell=x_i^n$, then $w_\ell^*=x_i^{-n}$ and the weight is unchanged. If $w_\ell\in \cB_j$, then $w_\ell^*\in \cA_j^\circ$ and hence $w_\ell^*$ is a linear combination of elements of $\cB_j$. Reversing the order of the syllables preserves the condition that neighboring syllables come from different free product factors. Therefore $w^*$ is a linear combination of reduced words with the same extent as $w$.
\end{proof}

\begin{lemma}\label{lem:span-products}
For every nonnegative integer $\vd$,
\[
 \algAs_{2\vd}
 \,=\,
  \spanop\{v^*u \ :\ u,v\in \algAs_\vd\}.
\]
\end{lemma}

\begin{proof}
The inclusion from right to left follows from Lemma \ref{lem:adjoint-invariant} and the subadditivity of extent under multiplication.

\smallskip

For the reverse inclusion, let $w=w_1\cdots w_m$ be a reduced word with $\ext(w)\leq 2\vd$. Split $w$ as $w=pq$ with $\ext(p)\leq \vd$ and $\ext(q)\leq \vd$. If the splitting point falls inside a syllable $x_i^n$, write $x_i^n=x_i^t x_i^{n-t}$, choosing $t$ with the same sign as $n$ so that the two new powers have the desired extents. No symbol from $\mathcal{B}_j$ needs to be split, since all weights are integral and such syllables have weight one. Consequently,  $p,q\in\algAs_\vd$ and
\[
   w=pq=(p^*)^*q.
\]
By Lemma \ref{lem:adjoint-invariant}, $p^*\in \algAs_\vd$, so $w$ belongs to the stated span.
\end{proof}

The following subgroup reduction is a standard consequence of
induction of representations.

\begin{lemma}\label{r:reduce-to-finite}
Let $\groupS$ be a subgroup of $\groupY$ and let $D$
be a unital $C^*$-algebra. For
$\alpha\in\{\min,\max\}$, the canonical map
\[
 D\otimes_\alpha C^*(\groupS)
 \longrightarrow D\otimes_\alpha C^*(\groupY)
\]
is injective. Consequently, for every
$A=A^*\in D\otimes\C[\groupS]$ and $\varepsilon\ge0$,
\[
 A\succeq\varepsilon I
 \quad\text{in }D\otimes_\alpha C^*(\groupS)
 \quad\Longleftrightarrow\quad
 A\succeq\varepsilon I
 \quad\text{in }D\otimes_\alpha C^*(\groupY).
\]
The same equivalence holds for strict positivity
in these $C^*$-algebras.

In particular, taking $D=M_K(C_u^*(\algA))$ and
$\alpha=\max$ recovers the corresponding equivalence
for pointwise strict positivity under all unital
$*$-representations of
$\algA\otimes\C[\groupS]$ and
$\algA\otimes\C[\groupY]$.
\end{lemma}

\begin{proof}
The inclusion $C^*(\groupS)\subseteq C^*(\groupY)$
is standard; see \cite[Proposition~2.5.8]{BO08}.
Injectivity of the minimal tensor product gives
the assertion for $\alpha=\min$.
For $\alpha=\max$, induction of representations
from $\groupS$ to $\groupY$, with $D$ acting
fiberwise, preserves the maximal tensor norm.
The positivity equivalences follow by spectral
permanence. The final assertion follows from
the characterization of strict positivity by
representations; see \cite[2.12.12]{Dixmier}.
\end{proof}

\section{Positive truncated kernels and representations}
\label{s:kernels-and-reps}

Let $k,d$ be positive integers, let $\groupS$ be a group and 
let $S\subseteq\groupS$.
 A linear map
\[
 \Theta:
 \algA_{2d}\otimes
 \operatorname{span}\{g:g\in\lfrac S\}
 \longrightarrow M_k(\C)
\]
is a \df{positive truncated kernel on $S$} if, for every
finite $G\subseteq S$,
\begin{equation}
\label{e:truncated-psd}
 \left[
  \Theta(p^*q\otimes u^{-1}v)
 \right]_{(p,u),(q,v)\in\cW_d\times G}
 \succeq0.
\end{equation}

We also write $\Theta(a,g)$ for $\Theta(a\otimes g)$,
and suppress the group factor when the group is trivial.
{If $\lly\subseteq \groupS$ is a monoid such that
$\groupS=\{u^{-1}v : u,v\in \lly\},$ then $\operatorname{span} \lfrac \lly = \C[\groupS].$ }

\begin{theorem}
\label{prop:main:rep}
Suppose $\groupS$ is an arbitrary discrete group and 
$\lly\subseteq\groupS$ is a monoid such that
\[
 \groupS=\{u^{-1}v:u,v\in\lly\}.
\]
If the linear map
\[
 \Theta:\algA_{2d}\otimes\C[\groupS]\longrightarrow M_k(\C)
\]
is a positive truncated kernel on $\lly$, then there
exist a Hilbert space $\HSE$, a unital $\ast$-representation
\[
 \pi:\algA\otimes\C[\groupS]\longrightarrow B(\HSE),
\]
and a bounded linear map $V:\C^k\to\HSE$ such that
\begin{equation}
\label{e:rep}
 \Theta(p^*q\otimes u^{-1}v)
 =V^*\pi\bigl((p^*q)\otimes(u^{-1}v)\bigr)V
\end{equation}
for all $p,q\in\algA_d$ and $u,v\in\groupS$.
Moreover, $V^*V=\Theta(\idA,\idY)$.
If $\groupS$ is at most countable, $\HSE$ can be chosen
to be separable.

More precisely, let $\cH_\Theta$ be the Hilbert space
obtained by quotienting
$\algA_d\otimes\C[\groupS]\otimes\C^k$ by the null
space of the positive semidefinite form
\[
 \langle p\otimes u\otimes\xi,
         q\otimes v\otimes\eta\rangle_\Theta
 =\langle\Theta(q^*p,v^{-1}u)\xi,\eta\rangle_{\C^k}
\]
and completing. Write $[p\otimes u\otimes\xi]$ for the
resulting classes, and let $\tau_\Theta$ be the canonical
unitary representation of $\groupS$ on $\cH_\Theta$,
\[
 \tau_\Theta(g)[p\otimes u\otimes\xi]
       =[p\otimes gu\otimes\xi].
\]
One can choose the realization in \eqref{e:rep} so that
\begin{equation}
\label{e:rep-canonical-action}
 \begin{gathered}
 \HSE=\cH_\Theta\otimes\ell^2(\N),
 \qquad
 V\xi=[\idA\otimes\idY\otimes\xi]\otimes e_0,\\
 \pi(\idA\otimes g)
       =\tau_\Theta(g)\otimes I_{\ell^2(\N)},
       \qquad g\in\groupS,
 \end{gathered}
\end{equation}
where $e_0$ is the first standard basis vector of
$\ell^2(\N)$. In this realization,
\begin{equation}
\label{e:rep-canonical-vectors}
 \pi(p\otimes g)V\xi
       =[p\otimes g\otimes\xi]\otimes e_0,
 \qquad p\in\algA_d,\quad g\in\groupS,\quad \xi\in\C^k.
\end{equation}

In the case that $\groupS=\{\idY\}$ is
trivial (so that
$\algA\otimes\C[\groupS]\cong\algA$), the realization
in \eqref{e:rep} can be chosen with $\HSE$ of dimension
at most
\[
\operatorname{lcm}
\bigl\{
1,2,\ldots,kN(d),
\dim\cA_1,\ldots,\dim\cA_s
\bigr\}.
\]
\end{theorem}

The proof of Theorem~\ref{prop:main:rep} uses the following two lemmas.
The first  lemma extends the positivity condition from finite
subsets of $\lly$ to finite subsets of $\groupS$.
For $\groupS=\Z$ and $\lly=\N_0$, it gives the passage
from unilateral to bilateral Toeplitz positivity.
It follows from  a standard fact from noncommutative
ring theory. See  \cite[Lemma~3.6]{KLM} for instance. A proof is given
here for the reader's convenience. 

\begin{lemma}
\label{l:inf-to-rep-plus}
Suppose that $\lly\subseteq\groupS$ is a monoid 
and $\groupS=\lfrac\lly$. If
\[
 \Theta:\algA_{2d}\otimes\C[\groupS]\longrightarrow M_k(\C)
\]
is a positive truncated kernel on $\lly$, then, for every
finite $G\subseteq\groupS$,
\begin{equation}
\label{e:Y-to-groupY}
 \left[
  \Theta(p^*q\otimes u^{-1}v)
 \right]_{(p,u),(q,v)\in\cW_d\times G}
 \succeq0.
\end{equation}
\end{lemma}

\begin{proof}
We first find $c\in\lly$ such that $cG\subseteq\lly$.
Write $G=\{g_1,\ldots,g_m\}$ and set $c_0=\idY$.
Inductively, write
\[
 c_{j-1}g_j=a_j^{-1}b_j,
 \qquad a_j,b_j\in\lly,
\]
and put $c_j=a_jc_{j-1}$.
Then $c_jg_i\in\lly$ for $1\le i\le j$.
Thus $c=c_m$ has the required property.

For $u,v\in G$,
\[
 (cu)^{-1}(cv)=u^{-1}v.
\]
Therefore the matrix in \eqref{e:Y-to-groupY} is the
positive matrix associated with the finite subset
$cG\subseteq\lly$, after reindexing $u$ as $cu$.
\end{proof}

\begin{lemma} %
\label{l:key-extend}
Let $\tau:G\to B(\cH)$ be a unitary representation of a
group $G$, and let $\cM,\cN$ be closed reducing subspaces
for $\tau$. Set
\[
 \widetilde{\cH}=\cH\otimes\ell^2(\N),
 \qquad
 \wtau(g)=\tau(g)\otimes I_{\ell^2(\N)},
\]
and identify $\cH$ with the first summand of
$\widetilde{\cH}$.
Write $\wtau_{\cM}$ and $\wtau_{\cN}$ for the restrictions
to $\widetilde{\cH}\ominus\cM$ and
$\widetilde{\cH}\ominus\cN$, respectively.

\begin{enumerate}[(i)]
\item There is a unitary
$V:\widetilde{\cH}\ominus\cM\to\widetilde{\cH}$ such that
\[
 V\wtau_{\cM}(g)=\wtau(g)V,\qquad g\in G.
\]

\item There is a unitary
$R:\widetilde{\cH}\ominus\cM
   \to\widetilde{\cH}\ominus\cN$ such that
\[
 R\wtau_{\cM}(g)=\wtau_{\cN}(g)R,\qquad g\in G.
\]

\item If $L:\cM\to\cN$ is unitary and intertwines
$\tau|_{\cM}$ and $\tau|_{\cN}$, then $L$ extends to a
unitary on $\widetilde{\cH}$ commuting with $\wtau$.

\item If $\cA$ is a finite-dimensional unital $C^*$-algebra
and $L:\cA\to B(\cM)$ is a unital $*$-representation
commuting with $\tau|_{\cM}$, then there is a unital
$*$-representation
\[
 \rho:\cA\to B(\widetilde{\cH}\ominus\cM)
\]
commuting with $\wtau_{\cM}$.
Consequently, $\sigma=L\oplus\rho$ is a unital
$*$-representation on $\widetilde{\cH}$ commuting with
$\wtau$.
\end{enumerate}
\end{lemma}

\begin{proof}
Note that $\cH=\{h\otimes e_0:h\in\cH\}$, where
$(e_n)$ is the standard orthonormal basis for $\ell^2(\N)$.
Since $\cM$ reduces $\tau$, we have the direct sum
decomposition
$\tau=\tau|_{\cM}\oplus\tau|_{\cM^\bot}$.
Consequently,
\[
\begin{split}
 \widetilde{\tau}|_{\widetilde{\cH}\ominus\cM}
 &\simu \tau|_{\cM^\bot}
       \oplus(\tau|_{\cM}\oplus\tau|_{\cM^\bot})
       \oplus(\tau|_{\cM}\oplus\tau|_{\cM^\bot})
       \oplus\cdots\\
 &\simu(\tau|_{\cM^\bot}\oplus\tau|_{\cM})
       \oplus(\tau|_{\cM^\bot}\oplus\tau|_{\cM})
       \oplus\cdots
 \simu\widetilde{\tau},
\end{split}
\]
where $\simu$ indicates unitary equivalence. Likewise,
$\widetilde{\tau}|_{\widetilde{\cH}\ominus\cN}
\simu\widetilde{\tau}$.
Hence there are unitary operators
\[
 V:\widetilde{\cH}\ominus\cM\longrightarrow\widetilde{\cH},
 \qquad
 R:\widetilde{\cH}\ominus\cM
       \longrightarrow\widetilde{\cH}\ominus\cN
\]
such that
\begin{equation}
 \label{e:Vtau}
 V\wtau_\cM(g)=\wtau(g)V,
 \qquad
 R\wtau_\cM(g)=\wtau_\cN(g)R,
 \qquad g\in G.
\end{equation}
This proves (i) and (ii).
For (iii), let $U=L\oplus R$ and observe that, for $g\in G$,
\[
\begin{split}
 U\wtau(g)
 &=(L\oplus R)
       (\tau(g)|_{\cM}\oplus\wtau_\cM(g))\\
 &=L\tau(g)|_{\cM}\oplus R\wtau_\cM(g)\\
 &=\tau(g)|_{\cN}L\oplus\wtau_\cN(g)R\\
 &=(\tau(g)|_{\cN}\oplus\wtau_\cN(g))(L\oplus R)
 =\wtau(g)U.
\end{split}
\]

To prove (iv), choose a unital $\ast$-representation
$\nu:\cA\to B(\ell^2(\N))$ and define
$\rho:\cA\to B(\widetilde{\cH}\ominus\cM)$ by
\[
 \rho(a)=V^*(I_\cH\otimes\nu(a))V,
 \qquad a\in\cA,
\]
 set $\widetilde{\nu}= I_\cH\otimes \nu$ 
 and observe, by equation~\eqref{e:Vtau} and the fact that $V$ is unitary, 
 \[
 \begin{split}
   \rho(a) \wtau_\cM(g)  &= V^* \widetilde{\nu}(a) V \wtau_\cM(g)
    = V^* \widetilde{\nu}(a) \wtau(g) V
  \\&    =V^* \wtau(g) \widetilde{\nu}(a) V
     = \wtau_{\cM} (g) (V^* \widetilde{\nu}(a) V)
      = \wtau_\cM(g) \rho(a).
  \end{split}
 \]
Define 
\[
 \sigma:\cA\longrightarrow
 B\bigl(\cM\oplus(\widetilde{\cH}\ominus\cM)\bigr)
 =B(\widetilde{\cH})
\]
by $\sigma(a)=L(a)\oplus\rho(a)$.
Thus $\sigma$ is a unital $\ast$-representation of $\cA$
on $\widetilde{\cH}$ and, for $a\in\cA$ and $g\in G$,
\begin{equation*}
\begin{split}
 \sigma(a)\widetilde{\tau}(g)
 &=[L(a)\oplus \rho(a)]
       [\tau(g)|_{\cM}\oplus\wtau_\cM(g)]\\[3pt]
 &=L(a)\tau(g)|_{\cM}
       \oplus  \rho(a)\wtau_\cM(g)  =  \tau(g)|_\cM L(a) \oplus \wtau_\cM(g) \rho(a)
  =   \wtau(g) \sigma(a).
 \end{split}
\end{equation*}
\end{proof}

\begin{proof}[Proof of Theorem~\ref{prop:main:rep}]
By Lemma~\ref{l:inf-to-rep-plus}, the kernel $\Theta:\algA_{2d}\otimes\C[\groupS]\to M_k(\C)$ is positive on every finite subset of $\groupS$.
 On $\algA_d\otimes\C[\groupS]\otimes\C^k,$ define
\[
 \langle p\otimes u\otimes\xi,q\otimes  v\otimes\eta\rangle_\Theta =\langle\Theta(q^*p,v^{-1}u)\xi,\eta\rangle_{\C^k}.
\]
By equation~\eqref{e:truncated-psd}, $\langle\cdot, \, \cdot\rangle_\Theta$  is a positive semidefinite sesquilinear form. Let
\[
 \cN_\Theta=\{\omega:\langle\omega,\omega\rangle_\Theta=0\},
\]
and let $\HSH_\Theta$ be the Hilbert-space completion of
\[
 (\algA_d\otimes\C[\groupS]\otimes\C^\kk)/\cN_\Theta
\]
 under the inner product induced by $\langle \cdot,\cdot\rangle_\Theta$
 and  denote the resulting equivalence classes by $[p\otimes u\otimes\xi].$ 
 If $\groupS$ is at most countable, then the finite-dimensionality of
$\algA_d$ and $\C^k$ implies that $\HSH_\Theta$ is separable.
The construction also applies when $\groupS$ is uncountable; in that
case no separability of $\HSH_\Theta$ is asserted.

\medskip

Define the canonical representation
$\tau=\tau_\Theta:\groupS\to B(\cH_\Theta)$ as follows. Since $(gv)^{-1}(gu)=v^{-1}u$ for $g, u, v\in \groupS,$
for $g\in \groupS$ setting
\[
 \tau(g)[p\otimes u\otimes\xi]=[p\otimes gu\otimes\xi]
\]
 and extending by linearity and continuity produces a unitary operator $\tau(g)\in B(\cH_\Theta).$
Since $\tau(gh)=\tau(g)\tau(h)$ as is readily checked,  $\tau$ is a unitary representation of $\groupS$  and extends to a unital 
$\ast$-representation of $\C[\groupS].$

\medskip

Fix $i\in\{1,\ldots,r\}$, and let $z=x_i$ or $z=x_i^{-1}$.
Before passing to the quotient, put
\[
 \cM_z^0
 =
 \operatorname{span}
 \{p\otimes u\otimes\xi:
   p\in\algA_{d-1}+z\algA_{d-1},\
   u\in\groupS,\ \xi\in\C^k\},
\]
and let $\cM_z$ be the closure of its image in $\cH_\Theta$.
Left multiplication by $z$  defines a bijection
\[
 L_z^0:\cM_{z^{-1}}^0\to\cM_z^0,
 \qquad
 L_z^0(p\otimes u\otimes\xi)=zp\otimes u\otimes\xi,
\]
with inverse $L_{z^{-1}}^0$.
For $p,q\in\algA_{d-1}+z^{-1}\algA_{d-1}$,
\[
\begin{aligned}
 &\langle L_z^0(p\otimes u\otimes\xi),
             L_z^0(q\otimes v\otimes\eta)\rangle_\Theta\\
 &\qquad=
 \langle\Theta(q^*z^*zp\otimes v^{-1}u)\xi,\eta\rangle\\
 &\qquad=
 \langle\Theta(q^*p\otimes v^{-1}u)\xi,\eta\rangle\\
 &\qquad=
 \langle p\otimes u\otimes\xi,
         q\otimes v\otimes\eta\rangle_\Theta.
\end{aligned}
\]
By sesquilinearity, this identity holds throughout
$\cM_{z^{-1}}^0$. Thus $L_z^0$ descends to the quotient
and extends to a unitary
\[
 L_z:\cM_{z^{-1}}\to\cM_z,
 \qquad L_z^{-1}=L_{z^{-1}}.
\]

\medskip

Fix $j\in\{1,\ldots,s\}.$ Before passing to the quotient, let
\[
 \cM_j^0=\operatorname{span}\{ap\otimes u\otimes\xi:a\in\mathcal A_j,\ p\in \algA_{d-1},\ u\in \groupS,\ \xi\in\C^\kk\}.
\]
Define $L_j^0:\cA_j\to\cL(\cM_j^0)$, where
$\cL(\cM_j^0)$ denotes the linear maps from $\cM_j^0$
to itself, by
\[
 L_j^0(a)\left(\sum_\ell b_\ell p_\ell\otimes u_\ell\otimes\xi_\ell\right) =\sum_\ell ab_\ell p_\ell\otimes u_\ell\otimes\xi_\ell.
\]
One has
\begin{equation*}
\langle L_j^0(a)\zeta,\omega\rangle_\Theta =\langle\zeta,L_j^0(a^*)\omega\rangle_\Theta. 
\end{equation*}
Put $C=\|a\|.$ Since $C^2 -a^*a\succeq0,$ there is $b\in\mathcal A_j$ such that $C^2-a^*a=b^*b.$ It follows that
\begin{equation*}
 \|L_j^0(a)\zeta\|_\Theta\le \|a\|\,\|\zeta\|_\Theta, \qquad \zeta\in \cM_j^0. 
\end{equation*}
Thus $L_j^0(a)$ descends to the quotient and extends to a bounded operator $L_j(a)$ on
\[
 \cM_j=\overline{\cM_j^0/(\cM_j^0\cap \cN_\Theta)}.
\]
The map $L_j:\mathcal A_j\to B(\cM_j)$ is a unital $\ast$-representation. 
 
 \medskip
 
We pause here to prove the finite-dimensional alternative when
$\groupS=\{\idY\}$. In this case, $\HSH_\Theta$ is a quotient of
\[
\algAs_d\otimes\C^k.
\]
Since $\dim\algAs_d=N(d)$, we have
\[
\dim\HSH_\Theta\leq kN(d).
\]
Set
\begin{equation}
\label{e:D}
H=kN(d)
\qquad\text{and}\qquad
D=
\operatorname{lcm}
\bigl\{
1,2,\ldots,H,
\dim\cA_1,\ldots,\dim\cA_s
\bigr\}.
\end{equation}
In particular,
\[
D\geq H\geq\dim\HSH_\Theta.
\]
Choose a Hilbert space $\cK_\Theta$ such that
\[
\cE_\Theta=\HSH_\Theta\oplus\cK_\Theta
\qquad\text{and}\qquad
\dim\cE_\Theta=D.
\]

For each $i\in\{1,\ldots,r\}$, the local left-multiplication
operator
\[
L_{x_i}:\cM_{x_i^{-1}}\longrightarrow\cM_{x_i}
\]
is unitary. Consequently,
\[
\dim\cM_{x_i^{-1}}=\dim\cM_{x_i},
\]
and hence
\[
\dim(\cE_\Theta\ominus\cM_{x_i^{-1}})
=
\dim(\cE_\Theta\ominus\cM_{x_i}).
\]
Thus  $L_{x_i}$ extends to a unitary
\[
U_i\in B(\cE_\Theta).
\]
Define
\[
\pi_{x_i}:\C[x_i,x_i^{-1}]\longrightarrow B(\cE_\Theta),
\qquad
\pi_{x_i}(x_i)=U_i.
\]

We next extend the local representations associated with the
finite-dimensional factors. Fix $j\in\{1,\ldots,s\}$ and put
\[
m_j=\dim\cM_j.
\]
Then
\[
0\leq m_j\leq\dim\HSH_\Theta\leq H.
\]

Suppose first that $m_j>0$. Since $m_j\leq H$, the definition of
$D$ implies that $m_j$ divides $D$. Hence
\[
D-m_j
=
\left(\frac{D}{m_j}-1\right)m_j.
\]
Recall that
\[
L_j:\cA_j\longrightarrow B(\cM_j)
\]
is a unital $*$-representation. Let
\[
\widehat\rho_j
=
L_j^{\oplus(D/m_j-1)}
\]
be the direct sum of $D/m_j-1$ copies of $L_j$. This is a unital
$*$-representation on a Hilbert space of dimension $D-m_j$.
After identifying that Hilbert space unitarily with
$\cE_\Theta\ominus\cM_j$, we obtain a unital $*$-representation
\[
\rho_j:
\cA_j\longrightarrow B(\cE_\Theta\ominus\cM_j).
\]
Define
\[
\pi_j(a)=L_j(a)\oplus\rho_j(a),
\qquad a\in\cA_j,
\]
with respect to the orthogonal decomposition
\[
\cE_\Theta
=
\cM_j\oplus(\cE_\Theta\ominus\cM_j).
\]

Suppose instead that $m_j=0$. Equip $\cA_j$ with the Hilbert-space
inner product
\[
\langle a,b\rangle_j=\tau_j(b^*a)
\]
and let
\[
\lambda_j:\cA_j\longrightarrow B(\cA_j)
\]
denote the left regular representation,
\[
\lambda_j(a)b=ab,
\qquad a,b\in\cA_j.
\]
The faithfulness of $\tau_j$ implies that this inner product is
positive definite, and the tracial property gives
\[
\langle\lambda_j(a)b,c\rangle_j
=
\langle b,\lambda_j(a^*)c\rangle_j.
\]
Thus $\lambda_j$ is a unital $*$-representation on a Hilbert space
of dimension $\dim\cA_j$. Since $\dim\cA_j$ divides $D$, the direct
sum
\[
\lambda_j^{\oplus(D/\dim\cA_j)}
\]
is a unital $*$-representation on a Hilbert space of dimension
$D$. After identifying this Hilbert space unitarily with
$\cE_\Theta$, let
\[
\pi_j:\cA_j\longrightarrow B(\cE_\Theta)
\]
denote the resulting representation.

By the universal property of the algebraic unital free product,
the representations
\[
\pi_{x_1},\ldots,\pi_{x_r},
\pi_1,\ldots,\pi_s
\]
determine a unital $*$-representation
\begin{equation}
\label{e:piTheta}
\pi_\Theta:\algA\longrightarrow B(\cE_\Theta).
\end{equation}
By construction, $\pi_\Theta$ extends each of the local
left-multiplication operators used above.
Define
\[
 V_\Theta:\C^k\to\cE_\Theta,
 \qquad
 V_\Theta\xi=[\idA\otimes\xi].
\]
Induction on extent, using the local left-multiplication
operators, gives
\[
 \pi_\Theta(q)V_\Theta\xi=[q\otimes\xi],
 \qquad q\in\algA_d.
\]
For inverse generators, use the fact that $U_i^*$ extends
$L_{x_i^{-1}}$, since $U_i$ extends the unitary $L_{x_i}$
between its initial and final subspaces.
Consequently, for $p,q\in\algA_d$,
\[
\begin{aligned}
 \langle V_\Theta^*\pi_\Theta(p^*q)V_\Theta\xi,\eta\rangle
 &=
 \langle[q\otimes\xi],[p\otimes\eta]\rangle_\Theta\\
 &=
 \langle\Theta(p^*q)\xi,\eta\rangle.
\end{aligned}
\]
Lemma~\ref{lem:span-products} now gives
\[
 \Theta(a)=V_\Theta^*\pi_\Theta(a)V_\Theta,
 \qquad a\in\algA_{2d}.
\]
Also $V_\Theta^*V_\Theta=\Theta(\idA)$.
This proves the finite-dimensional alternative, including
the stated dimension bound. 
 
\medskip

We now return to the general case and construct the canonical amplified realization for an
arbitrary group $\groupS$, including the trivial group. Set
\[
 \HSE=\widetilde{\cH}=\cH_\Theta\otimes\ell^2(\N),
 \qquad
 \wtau(g)=\tau(g)\otimes I_{\ell^2(\N)},
\]
and identify $\cH_\Theta$ with its first summand, by
$h\mapsto h\otimes e_0$.
Since $\tau(g)$ changes only the group coordinate, it
maps each of $\cM_z$ and $\cM_{z^{-1}}$ onto itself.
These spaces therefore reduce $\tau(\groupS)$.
For $p\in\algA_{d-1}+z^{-1}\algA_{d-1}$,
$u,g\in\groupS$ and $\xi\in\C^k$,
\[
\begin{aligned}
 L_z\tau(g)[p\otimes u\otimes\xi]
 &=[zp\otimes gu\otimes\xi]\\
 &=\tau(g)L_z[p\otimes u\otimes\xi].
\end{aligned}
\]
By linearity and continuity, this proves the intertwining
identity
\[
 L_z\tau(g)|_{\cM_{z^{-1}}}
       =\tau(g)|_{\cM_z}L_z,
 \qquad g\in\groupS.
\]
Hence, by Lemma~\ref{l:key-extend}(iii), for each
$i\in\{1,\ldots,r\}$ there is a unitary operator $U_i$
on $\widetilde{\cH}$ extending $L_{x_i}$ and satisfying
\begin{equation}
\label{e:[U,tau]=0}
 U_i\wtau(g)=\wtau(g)U_i,
 \qquad g\in\groupS.
\end{equation}
Since $U_i$ extends the unitary
$L_{x_i}:\cM_{x_i^{-1}}\to\cM_{x_i}$, its adjoint
$U_i^*$ extends $L_{x_i^{-1}}$.
As $F_r$ is the free group on $x_1,\ldots,x_r$, the
unitaries $U_1,\ldots,U_r$ determine a unitary
representation of $F_r$ that commutes with $\wtau$.

Similarly, $\tau(g)$ maps $\cM_j$ onto itself, so
$\cM_j$ reduces $\tau$. For $a,b\in\cA_j$,
$p\in\algA_{d-1}$, $u,g\in\groupS$ and $\xi\in\C^k$,
\[
\begin{aligned}
 L_j(a)\tau(g)[bp\otimes u\otimes\xi]
 &=[abp\otimes gu\otimes\xi]\\
 &=\tau(g)L_j(a)[bp\otimes u\otimes\xi].
\end{aligned}
\]
Again, linearity and continuity give
\[
 L_j(a)\tau(g)|_{\cM_j}
       =\tau(g)|_{\cM_j}L_j(a),
 \qquad a\in\cA_j,\quad g\in\groupS.
\]
By Lemma~\ref{l:key-extend}(iv), there is a unital
$\ast$-representation
\[
 \rho_j:\cA_j\longrightarrow
       B(\widetilde{\cH}\ominus\cM_j)
\]
such that the unital $\ast$-representation
$\sigma_j=L_j\oplus\rho_j$ commutes with $\wtau$:
\begin{equation}
\label{e:[sigma,tau]=0}
 \sigma_j(a)\wtau(g)=\wtau(g)\sigma_j(a),
 \qquad a\in\cA_j,\quad g\in\groupS.
\end{equation}

\medskip

The universal property of the algebraic unital free
product gives a unital $\ast$-representation
\[
 \sigma:\algA\longrightarrow B(\widetilde{\cH})
\]
extending the representations determined by the $U_i$
and the $\sigma_j$. It commutes with $\widetilde\tau$
by \eqref{e:[U,tau]=0} and \eqref{e:[sigma,tau]=0}.
Consequently,
\begin{equation}
\label{e:pi=sigma-tau}
 \pi(a\otimes g)=\sigma(a)\widetilde\tau(g)
\end{equation}
defines a unital $\ast$-representation of
$\algA\otimes\C[\groupS]$, with
$\pi(\idA\otimes g)=\tau_\Theta(g)\otimes I_{\ell^2(\N)}$.
Define
\[
 V:\C^k\longrightarrow\widetilde{\cH},
 \qquad
 V\xi=[\idA\otimes\idY\otimes\xi]\otimes e_0.
\]
This map is bounded, since
\[
 \|V\xi\|^2
 =\langle\Theta(\idA,\idY)\xi,\xi\rangle
 \le\|\Theta(\idA,\idY)\|\,\|\xi\|^2.
\]
Thus the realization satisfies \eqref{e:rep-canonical-action}.

We verify the vector identity
\eqref{e:rep-canonical-vectors}. In the following
calculation the kernel classes are viewed in the first
summand of $\widetilde{\cH}$. For the identity word,
\[
 \widetilde\tau(u)V\xi=[\idA\otimes u\otimes\xi].
\]
Suppose that the identity \eqref{e:rep-canonical-vectors} has been established for words
of extent less than $m\le d$. A nonempty reduced word
$p$ of extent $m$ can be written as $p=zq$, where
$z=x_i$ or $x_i^{-1}$ and $q\in\cW_{m-1}$, or as
$p=aq$, where $a\in\cB_j\subseteq\cA_j$ and
$q\in\cW_{m-1}$. In the first case,
$[q\otimes u\otimes\xi]\in\cM_{z^{-1}}$, and
$\sigma(z)$ extends $L_z$. In the second case,
$[q\otimes u\otimes\xi]\in\cM_j$, and
$\sigma(a)$ extends $L_j(a)$. Thus, respectively,
\[
\begin{aligned}
 \sigma(zq)\widetilde\tau(u)V\xi
 &=\sigma(z)[q\otimes u\otimes\xi]
   =[zq\otimes u\otimes\xi],\\
 \sigma(aq)\widetilde\tau(u)V\xi
 &=\sigma(a)[q\otimes u\otimes\xi]
   =[aq\otimes u\otimes\xi],
\end{aligned}
\]
completing the induction. By linearity and equation~\eqref{e:pi=sigma-tau},
\[
 \pi(p\otimes u)V\xi
 =\sigma(p)\widetilde\tau(u)V\xi
 =[p\otimes u\otimes\xi],
 \qquad p\in\algA_d,\quad u\in\groupS,
\]
which is \eqref{e:rep-canonical-vectors} under the
first-summand identification.

Finally, for $p,q\in\algA_d$, $u,v\in\groupS$ and
$\xi,\eta\in\C^k$,
\[
\begin{aligned}
 &\langle V^*\pi\bigl((p^*q)\otimes(u^{-1}v)\bigr)
                V\xi,\eta\rangle\\
 &\qquad=\langle\sigma(q)\widetilde\tau(v)V\xi,
                    \sigma(p)\widetilde\tau(u)V\eta\rangle\\
 &\qquad=\langle[q\otimes v\otimes\xi],
                    [p\otimes u\otimes\eta]\rangle_\Theta\\
 &\qquad=\langle\Theta(p^*q,u^{-1}v)\xi,\eta\rangle,
\end{aligned}
\]
proving  \eqref{e:rep} for all $u,v\in\groupS$.
Taking $p=q=\idA$ and $u=v=\idY$ gives
$V^*V=\Theta(\idA,\idY)$.
If $\groupS$ is at most countable, then $\cH_\Theta$
and its countable amplification $\HSE$ are separable.
\end{proof}

\section{Proofs of the main results}
\label{s:proofs}
In this section we prove Theorem~\ref{thm:main} and
Corollaries~\ref{cor:group-version} and \ref{c:main+finiteY}.
The proofs of Theorem~\ref{t:main+maxY}
and Corollary~\ref{t:main+amenableY} are given in
Section~\ref{s:operator-coefficients}.

\subsection{Proof of Theorem~\ref{thm:main}}

Enumerate $\cW_\vd$, 
\[
 \cW_\vd 
 \,=\,
 \{q_1,\ldots,q_{N(\vd)}\}.
\]
Let $\algA_{2\vd}^{*}$ be the set of linear functionals on $\algA_{2\vd}.$   Define
\[
 \cO_\vd
 \,=\,
 \left\{T^\lambda=\big(\lambda(q_\alpha^*q_\beta)\big)_{\alpha,\beta=1}^{N(\vd)} \ :\ \lambda\in \algAs_{2\vd}^*\right\}
\,\subseteq\,  M_{N(\vd)}(\C).
\]
By Lemma \ref{lem:span-products}, the map $\lambda\mapsto T^\lambda$ is injective. The space $\cO_\vd$ is an operator system. Indeed, it is self-adjoint because
\[
   (T^\lambda)^*=T^{\lambda^\sharp},
   \qquad
   \lambda^\sharp(a)=\overline{\lambda(a^*)},
\]
and it contains the identity because Lemma \ref{lem:orthonormal-reduced-words} gives
\[
   T^\tau=I_{N(\vd)}.
\]

 Let
\[
   A
    \,=\,
    \sum_{\ell=1}^{N(2\vd)} A_\ell\otimes r_\ell
   \in B(\cH)\otimes\algAs_{2 \vd}
\]
be positive on finite-dimensional representations. Using the fact that the map $\algA_{2 \vd}^* \ni \lambda \mapsto T^{\lambda}$ is one-to-one,
{we  obtain a well-defined map}
$  \Phi_A:\cO_\vd\to B(\cH) $
{defined} by
\[
   \Phi_A(T^\lambda)=(\id_{B(\cH)}\otimes\lambda)(A)
   =\sum_{\ell=1}^{N(2\vd)} A_\ell\lambda(r_\ell).
\]

\begin{lemma}\label{lem:Phi-cp}
The map $\Phi_A:\cO_\vd\to B(\cH)$ is completely positive.
\end{lemma}

\begin{proof}
Let
\[
   X=[T^{\lambda_{\mu,\nu}}]_{\mu,\nu=1}^k\in M_k(\cO_\vd)
\]
be positive. Define
\[
   \Theta:\algAs_{2\vd}\to M_k(\C),
   \qquad
   \Theta(a)=\big(\lambda_{\mu,\nu}(a)\big)_{\mu,\nu=1}^k.
\]
The positivity of $X$ is precisely the assertion that $\Theta$ is a positive truncated kernel. Let $\cE_\Theta$, $\pi_\Theta$ and $V$
be as in the finite-dimensional alternative in Theorem~\ref{prop:main:rep}. In particular, 
\[
   \Theta(a)=V^*\pi_\Theta(a)V,
   \qquad a\in\algAs_{2 \vd}.
\]
Under the standard identification $M_k(B(\cH))\cong B(\cH\otimes\C^k)$,
\[
\begin{aligned}
   \Phi_A^{(k)}(X)
   &=\sum_{\ell=1}^{N(2\vd)} A_\ell\otimes \Theta(r_\ell)\\
   &=(I_\cH\otimes V)^*
     \left(\sum_{\ell=1}^{N(2\vd)} A_\ell\otimes \pi_\Theta(r_\ell)\right)
     (I_\cH\otimes V)\\
   &=(I_\cH\otimes V)^*(\id_{B(\cH)}\otimes\pi_\Theta)(A)(I_\cH\otimes V).
\end{aligned}
\]
The representation $\pi_\Theta$ is finite-dimensional. Hence the assumed positivity of $A$ implies
\[
   (\id_{B(\cH)}\otimes\pi_\Theta)(A)\succeq 0.
\]
Therefore $\Phi_A^{(k)}(X)\succeq 0$. Since $k$ was arbitrary, $\Phi_A$ is completely positive.
\end{proof}

By Lemma \ref{lem:Phi-cp}, $\Phi_A$ is completely positive on the operator system $\cO_\vd\subseteq M_{N(\vd)}(\C)$. By Arveson's Extension Theorem \cite[Theorem 7.5]{Pa}, $\Phi_A$ extends to a completely positive map
\[
   \widetilde\Phi_A:M_{N(\vd)}(\C)\to B(\cH).
\]
Let $E_{\alpha,\beta}$ denote the matrix units  in $M_{N(\vd)}(\C)$. The Choi matrix
\[
   C=\big(\widetilde\Phi_A(E_{\alpha,\beta})\big)_{\alpha,\beta=1}^{N(\vd)}
   \in M_{N(\vd)}(B(\cH))
\]
is positive. As a result, there exist operators
\[
   B_1,\ldots,B_{N(\vd)}\in B(\cH,\cE),
   \qquad \cE=\cH\otimes\C^{N(\vd)},
\]
such that
\[
   \widetilde\Phi_A(E_{\alpha,\beta})=B_\alpha^*B_\beta.
\]
Define
\[
   B=\sum_{\alpha=1}^{N(\vd)}B_\alpha\otimes q_\alpha
   \in B(\cH,\cE)\otimes\algAs_\vd.
\]
Then
\[
   B^*B=
   \sum_{\alpha,\beta=1}^{N(\vd)}B_\alpha^*B_\beta\otimes q_\alpha^*q_\beta.
\]
Let $\lambda\in\algAs_{2\vd}^*$. Applying $\id\otimes\lambda$ gives
\[
\begin{aligned}
   (\id\otimes\lambda)(B^*B)
   &=\sum_{\alpha,\beta=1}^{N(\vd)}B_\alpha^*B_\beta\lambda(q_\alpha^*q_\beta)\\
   &=\sum_{\alpha,\beta=1}^{N(\vd)}\widetilde\Phi_A(E_{\alpha,\beta})\lambda(q_\alpha^*q_\beta)\\
   &=\widetilde\Phi_A(T^\lambda)\\
   &=\Phi_A(T^\lambda)\\
   &=(\id\otimes\lambda)(A).
\end{aligned}
\]
Since the linear functionals on $\algAs_{2\vd}$ separate points, the maps $\id\otimes\lambda$ separate points of the algebraic tensor product $B(\cH)\otimes\algAs_{2\vd}$. Hence $A=B^*B$.

To establish the  bound given in Remark~\ref{r:rep-bound}, in the case that $\cH$ has dimension $K,$ it suffices to check that $\Phi_A$
is $K$-positive \cite[Theorem~6.1]{Pa}. In that case, take $k=K$; the representation $\pi_\Theta$
is on a space of dimension at most $D$ from
equation~\eqref{e:D}.

\subsection{Proof of Corollary~\ref{cor:group-version}}
For each finite group $\mathbb G_j$, take
	\[
	\mathcal A_j=\mathbb C[\mathbb G_j]
	\]
	and equip $\mathbb C[\mathbb G_j]$ with its canonical trace 
	\[
	\tau_j(g)=
	\begin{cases}
		1, & g=\idGroupAlg,\\
		0, & g\neq \idGroupAlg.
	\end{cases}
	\]
	Then
	\[
	\mathbb C[\mathbb G_j]
	=
	\mathbb C\idA\oplus \mathbb C[\mathbb G_j]^\circ,
	\qquad
	\mathbb C[\mathbb G_j]^\circ=\ker\tau_j,
	\]
	and the set
	\[
	\mathbb G_j\setminus\{ \idGroupAlg\}
	\]
	is an orthonormal basis of $\mathbb C[\mathbb G_j]^\circ$ with respect to the inner
	product
	\[
	\langle a,b\rangle_j=\tau_j(b^*a).
	\]
	Moreover,
	\[
	\mathbb C[\mathbb G]
	\, \cong\, 
	\mathbb C[F_r]*\mathbb C[\mathbb G_1]*\cdots *\mathbb C[\mathbb G_s].
	\]
	A unitary representation of $\mathbb G$ is the same thing as a unital
	$*$-representation of the group algebra $\mathbb C[\mathbb G]$.  Hence the
	positivity assumption in the corollary is precisely the positivity assumption
	appearing in Theorem~\ref{thm:main}.  
With these choices of bases, $\algA_\vd$ is precisely the span of the group elements of extent at most $\vd$.
	Applying Theorem~\ref{thm:main} gives a
	Hilbert space $\mathcal E$ and an element
	\[
	B\in B(\mathcal H,\mathcal E)\otimes \algA_\vd
	\]
	such that
	\[
	A=B^*B.
	\]

\subsection{Proof of Corollary~\ref{c:main+finiteY}}
Since $\groupS$ is finite, its group algebra is a finite-dimensional $C^*$-algebra. Writing $\widehat{\groupS}$ for the set of irreducible unitary representations $\rho$ of $\groupS$, and $n_\rho=\dim\rho$, the Fourier transform gives a $*$-isomorphism
\[
\mathbb C[\groupS]
\cong
\bigoplus_{\rho\in\widehat{\groupS}} M_{n_\rho}(\mathbb C),
\qquad
y\longmapsto \bigoplus_{\rho\in\widehat{\groupS}}\rho(y).
\]

Write
\[
A
=
\sum_{r\in \cW_{2\vd}}\sum_{y\in \groupS}
A_{r,y}\otimes r\otimes y.
\]
For each $\rho\in\widehat{\groupS}$, define
\[
A_\rho
=
\sum_{r\in \cW_{2\vd}}
\left(
\sum_{y\in \groupS} A_{r,y}\otimes \rho(y)
\right)\otimes r
\in
B(\cH)\otimes M_{n_\rho}(\mathbb C)\otimes\algAs_{2\vd}.
\]

Let
$
\sigma:\algA\to B(\cF)
$
be a finite-dimensional unital $*$-representation. Then
\[
\pi_{\sigma,\rho}(a\otimes y)
=
\sigma(a)\otimes \rho(y)
\]
is a finite-dimensional unital $*$-representation of
$\algA\otimes\mathbb C[\groupS]$.
After interchanging the tensor factors $\cF$ and
$\C^{n_\rho}$, the assumed positivity of $A$ gives
\[
(\operatorname{id}\otimes \sigma)(A_\rho)\succeq 0.
\]
That is, $A_\rho$ satisfies the hypotheses of Theorem~\ref{thm:main}.
Put $\cE=\cH\otimes\C^{N(\vd)}$.
Applying Theorem~\ref{thm:main} with coefficient Hilbert
space $\cH\otimes\C^{n_\rho}$, and identifying
\[
 (\cH\otimes\C^{n_\rho})\otimes\C^{N(\vd)}
 \cong \cE\otimes\C^{n_\rho},
\]
we obtain
\[
 B_\rho
 =
 \sum_{p\in\cW_\vd}B_{\rho,p}\otimes p
 \in
 B(\cH\otimes\C^{n_\rho},\cE\otimes\C^{n_\rho})
 \otimes\algA_\vd
\]
such that $A_\rho=B_\rho^*B_\rho$.

The Fourier transform induces a linear isomorphism
\[
 B(\cH,\cE)\otimes\C[\groupS]
 \longrightarrow
 \bigoplus_{\rho\in\widehat{\groupS}}
 B(\cH\otimes\C^{n_\rho},\cE\otimes\C^{n_\rho}).
\]
Consequently, for each $p\in\cW_\vd$, there are unique
operators $B_{p,y}\in B(\cH,\cE)$, indexed by
$y\in\groupS$, such that
\[
 B_{\rho,p}
 =
 \sum_{y\in\groupS}B_{p,y}\otimes\rho(y),
 \qquad \rho\in\widehat{\groupS}.
\]
Define
\[
 B
 =
 \sum_{p\in\cW_\vd}\sum_{y\in\groupS}
 B_{p,y}\otimes p\otimes y
 \in B(\cH,\cE)\otimes\algA_\vd\otimes\C[\groupS].
\]

With respect to the Fourier transform, the
$\rho$-block of $B^*B$ is
\[
B_\rho^*B_\rho=A_\rho,
\]
which is exactly the $\rho$-block of $A$. Hence
\begin{equation}\label{eq:A=B^*B}
A=B^*B.
\end{equation}

Finally, suppose $\dim\cH=K$, and set
\[
 D_{\groupS}
 =
 \operatorname{lcm}
 \bigl\{
 1,2,\ldots,K|\groupS|N(\vd),
 \dim\cA_1,\ldots,\dim\cA_s
 \bigr\}.
\]
For $\rho\in\widehat{\groupS}$, Remark~\ref{r:rep-bound}
shows that it suffices to test $A_\rho$ on representations
of $\algA$ of dimension at most
\[
 D_\rho
 =
 \operatorname{lcm}
 \bigl\{
 1,2,\ldots,Kn_\rho N(\vd),
 \dim\cA_1,\ldots,\dim\cA_s
 \bigr\}.
\]
Since $\vd\ge1$, we have $\dim\cA_j\le N(\vd)$, and hence
\[
 D_\rho=\operatorname{lcm}\{1,\ldots,Kn_\rho N(\vd)\}.
\]
Also $n_\rho^2\le|\groupS|$, so
\[
 n_\rho D_\rho
 =
 \operatorname{lcm}
 \{n_\rho,2n_\rho,\ldots,Kn_\rho^2N(\vd)\}
 \quad\text{divides}\quad D_{\groupS}.
\]
The representations $\pi_{\sigma,\rho}$ needed in the
preceding argument therefore have dimension at most
$D_{\groupS}$. This proves the asserted testing bound.

\section{Proofs for tensor products}
\label{s:operator-coefficients}
We prove Theorem~\ref{t:main+maxY} and deduce
Corollary~\ref{t:main+amenableY}. We also provide an example delineating 
the distinction between spatial and maximal tensor-product positivity.

\subsection{Proof of Theorem~\ref{t:main+maxY}}

\begin{proof}
By \eqref{e:maxY-positive}, 
since $ B(\cH)\otm C_u^*(\algA)\otm C^*(\groupS)$ is a 
$C^*$-algebra,
there exists
$\varepsilon>0$ such that $A\succeq\varepsilon I$.

We use the separation and GNS argument of
\cite[Propositions~3.3 and~4.1]{HM04}, with the truncated
representation construction of Theorem~\ref{prop:main:rep}.
Put
\[
 \mathcal V_j=B(\cH)\otimes\algA_j\otimes\C[\groupS]
 \quad(j=\vd,2\vd),
 \qquad
 \Sigma_\vd=\left\{\sum_{k=1}^m F_k^*F_k:
                     F_k\in\mathcal V_\vd,\ m\ge1\right\},
\]
where tensor products without subscripts are algebraic. 

First, $I$ is an algebraic order unit for $\Sigma_\vd$ in
$(\mathcal V_{2\vd})_{\rm sa}$, i.e.,
\begin{equation}\label{e:maxY-order-unit}
 \text{for every }Z=Z^*\in\mathcal V_{2\vd}
 \text{ there is }R>0\text{ with }RI\pm Z\in\Sigma_\vd.
\end{equation}
 To see this, let
$F=b\otimes w\otimes g$, with $w=w_1\cdots w_m\in\cW_\vd$
a reduced word. Put $c_i=\|w_i\|^2$ and
$t_i=(c_i\idA-w_i^*w_i)^{1/2}$. Each $t_i$ belongs to
the same finite-dimensional factor as $w_i$, or is zero
when $w_i$ belongs to a free-group factor. The telescoping
identity
\[
 \left(\prod_i c_i\right)\idA-w^*w
 =\sum_i\left(\prod_{h<i}c_h\right)
       (t_iw_{i+1}\cdots w_m)^*(t_iw_{i+1}\cdots w_m)
\]
shows that the difference on the left is a sum of squares
of elements of $\algA_\vd$. Since $g$ is unitary and
$\|b\|^2I-b^*b$ is a square in $B(\cH)$, it follows
that $RI-F^*F\in\Sigma_\vd$ for some $R>0$.
The identity
\[
 m\sum_iF_i^*F_i-\left(\sum_iF_i\right)^*
                       \left(\sum_iF_i\right)
 =\sum_{i<j}(F_i-F_j)^*(F_i-F_j)
\]
extends this assertion to every $F\in\mathcal V_\vd$.
Polarization and Lemma~\ref{lem:span-products} now give
\eqref{e:maxY-order-unit}.

Fix $0<\delta<\varepsilon$ and suppose that
$A-\delta I\notin\Sigma_\vd$. By algebraic order-unit
separation, as in \cite[Proposition~3.3]{HM04}, there is
a complex-linear functional $\Lambda$ on $\mathcal V_{2\vd}$
such that
\begin{equation}\label{e:maxY-separation}
 \Lambda(F^*F)\ge0\quad(F\in\mathcal V_\vd),
 \qquad \Lambda(I)=1,
 \qquad \Lambda(A-\delta I)\le0.
\end{equation}

Equip $\mathcal V_\vd$ with the positive semidefinite form
$\langle F,G\rangle_\Lambda=\Lambda(G^*F)$, and let
$\cK_0$ be the Hilbert-space completion of the quotient
by its null space. Denote the classes by $[F]$ and put
$\xi=[I]$. Left multiplication in the two outer factors
gives commuting representations
\[
 \kappa_0(b)[F]=[(b\otimes\idA\otimes\idY)F],
 \qquad
 \upsilon_0(g)[F]=[(I_\cH\otimes\idA\otimes g)F].
\]
Here $\kappa_0$ is a unital $*$-representation of $B(\cH)$:
its boundedness follows from the square root of
$\|b\|^2I-b^*b$. Also $\upsilon_0$ is a unitary representation
of $\groupS$. Set
\[
 G=\mathcal U(B(\cH))\times\groupS,
 \qquad \rho(u,g)=\kappa_0(u)\upsilon_0(g),
\] 
where $\mathcal U(B(\cH))$ is the group of unitary operators
on $\cH$. Regard $G$ as a discrete group. Since the
unitaries linearly span $B(\cH)$, the canonical unital
$*$-homomorphism
\[
 \C[G]\longrightarrow B(\cH)\otimes\C[\groupS],
 \qquad (u,g)\longmapsto u\otimes g,
\]
is surjective.

To apply Theorem~\ref{prop:main:rep}, take $k=1$,
take the monoid to be $G$ itself, and define the
complex-linear map
\[
 \Theta:\algA_{2\vd}\otimes\C[G]\longrightarrow\C,
 \qquad
 \Theta(a,(u,g))=\Lambda(u\otimes a\otimes g).
\]
The positivity of $\Lambda$ implies that $\Theta$ is
a positive truncated kernel. Under the canonical
identification of $\cH_\Theta$ with $\cK_0$, the
representation $\tau_\Theta$ is $\rho$.

Put
\[
 \widetilde{\cK}=\cK_0\otimes\ell^2(\N),
 \qquad
 \widetilde\rho(s)=\rho(s)\otimes I_{\ell^2(\N)},
 \qquad s\in G.
\]
Identify $\cK_0$ with the first summand of
$\widetilde{\cK}$ and continue to write $\xi=[I]$ for
the vector $[I]\otimes e_0$.
The canonical realization in
Theorem~\ref{prop:main:rep} provides a unital
$*$-representation
\[
 \pi:\algA\otimes\C[G]\longrightarrow B(\widetilde{\cK})
\]
such that
\[
 \Theta(p^*q,s^{-1}t)
 =\bigl\langle
       \pi\bigl((p^*q)\otimes(s^{-1}t)\bigr)\xi,\xi
   \bigr\rangle,
 \qquad p,q\in\algA_\vd,\quad s,t\in G,
\]
and
\[
 \pi(\idA\otimes s)=\widetilde\rho(s),
 \qquad s\in G.
\]

Let $e_G=(I_\cH,\idY)$ be the identity of $G$, and define
\[
 \sigma:\algA\longrightarrow B(\widetilde{\cK}),
 \qquad
 \sigma(a)=\pi(a\otimes e_G).
\]
Since the two tensor factors commute in
$\algA\otimes\C[G]$, we have
\[
 \sigma(a)\widetilde\rho(u,g)
       =\widetilde\rho(u,g)\sigma(a),
 \qquad a\in\algA,\quad (u,g)\in G.
\]
Furthermore, \eqref{e:rep-canonical-vectors} gives
\[
 \sigma(q)\widetilde\rho(u,g)\xi
       =[u\otimes q\otimes g],
 \qquad q\in\algA_\vd,\quad (u,g)\in G,
\]
where the class on the right is viewed in the first
summand of $\widetilde{\cK}$.

Set
\[
 \kappa(b)=\kappa_0(b)\otimes I_{\ell^2(\N)},
 \qquad
 \upsilon(g)=\upsilon_0(g)\otimes I_{\ell^2(\N)},
 \qquad b\in B(\cH),\quad g\in\groupS.
\]
These representations have commuting ranges, and
\[
 \widetilde\rho(u,g)
 =\kappa(u)\upsilon(g)
 =\upsilon(g)\kappa(u),
 \qquad (u,g)\in G.
\]
Since $\sigma$ commutes with $\widetilde\rho$, taking
$(u,\idY)$ and $(I_\cH,g)$ shows that $\sigma(a)$
commutes with every $\kappa(u)$, $u$ unitary, and every
$\upsilon(g)$. The unitaries linearly span $B(\cH)$,
so $\sigma(\algA)$ commutes with
$\kappa(B(\cH))$ and $\upsilon(\groupS)$.

The representation $\sigma$ extends to $C_u^*(\algA)$,
and $\upsilon$ extends to $C^*(\groupS)$; their ranges
still commute with each other and with the range of
$\kappa$. By the universal property of the maximal
tensor product,
\[
 \Pi(b\otimes a\otimes g)
       =\kappa(b)\sigma(a)\upsilon(g)
\]
therefore extends to a unital $*$-representation of
$B(\cH)\otm C_u^*(\algA)\otm C^*(\groupS)$ on
$\widetilde{\cK}$.

By linearity in the unitary coefficient $u$, the vector
identity supplied by the theorem becomes
\[
 \sigma(q)\kappa(b)\upsilon(g)\xi
      =[b\otimes q\otimes g],
 \qquad b\in B(\cH),\quad q\in\algA_\vd,
        \quad g\in\groupS.
\]
In particular, for $b,c\in B(\cH)$,
$p,q\in\algA_\vd$ and $g,h\in\groupS$,
\[
\begin{aligned}
 &\bigl\langle
       \Pi(b^*c\otimes p^*q\otimes g^{-1}h)\xi,\xi
       \bigr\rangle\\
 &\qquad=\bigl\langle
       \kappa(c)\sigma(q)\upsilon(h)\xi,
       \kappa(b)\sigma(p)\upsilon(g)\xi\bigr\rangle\\
 &\qquad=\bigl\langle
       [c\otimes q\otimes h],
       [b\otimes p\otimes g]\bigr\rangle_\Lambda\\
 &\qquad=\Lambda(b^*c\otimes p^*q\otimes g^{-1}h).
\end{aligned}
\]
Taking $b=I_\cH$ and $g=\idY$, and using
Lemma~\ref{lem:span-products}, shows that  {the elements $b^*c\otimes p^*q\otimes g^{-1}h$}
span $\mathcal V_{2\vd}$. Hence
\[
 \Lambda(Z)=\langle\Pi(Z)\xi,\xi\rangle,
 \qquad Z\in\mathcal V_{2\vd}.
\]
Since $\|\xi\|^2=\Lambda(I)=1$, the inequality
$A\succeq\varepsilon I$ in the maximal tensor product
implies
\[
 \Lambda(A-\delta I)
 =\langle\Pi(A-\delta I)\xi,\xi\rangle
 \ge(\varepsilon-\delta)\|\xi\|^2
 =\varepsilon-\delta>0,
\]
contrary to \eqref{e:maxY-separation}. Therefore
$A-\delta I\in\Sigma_\vd$.

Because $A-\delta I\in\Sigma_\vd,$ 
\[
 A-\delta I=\sum_{j=1}^m F_j^*F_j,
 \qquad F_j\in\mathcal V_\vd.
\]
Stacking the $F_j$ in a column produces a finite
polynomial with coefficients in $B(\cH,\cH^{\oplus m})$
whose square is $A-\delta I$. Appending the constant row
$\sqrt\delta I_\cH$ gives
\[
 B=\begin{bmatrix}
      F_1\\ \vdots\\ F_m\\ \sqrt\delta I
    \end{bmatrix},
 \qquad B^*B=A.
\]
Every row has $\algA$-extent at most $\vd$, and their
group supports have finite union. Thus, with $n=m+1$,
this is a factorization of the asserted form.

For the final assertion, let {$\mathcal{F}$}  be the finite group
support of the factor {$B$}. The common-left-denominator
argument from Lemma~\ref{l:inf-to-rep-plus} gives
$c\in\lly$ with $c\cF\subseteq\lly$. Replacing every
group label $g$ by $cg$ leaves the square unchanged,
because $(cg)^{-1}(ch)=g^{-1}h$. Enlarging $M'$ if
necessary gives {$c\mathcal{F} \subseteq\lly_{\le M'}.$} 
\end{proof}

\subsection{Proof of Corollary~\ref{t:main+amenableY}}
{It turns out, building  on work of Kirchberg \cite{Kir93, Kir94}, results of
Pisier \cite[Theorems~0.1 and~0.2]{Pi96} say that the min and max
tensor products agree  in the context of Corollary~\ref{t:main+amenableY}
allowing for  an application of Theorem~\ref{t:main+maxY}. The Hilbert spaces 
in \cite{Pi96} are assumed separable,  a restriction 
that Lemma~\ref{lem:min-max-arbitrary-H} below says can be relaxed in our
setting.}

\begin{lemma}
\label{lem:min-max-arbitrary-H}
If  $C$ is  a unital $C^*$-algebra such that
\[
 B(\ell^2) \otimes_{\min}C
 =
 B(\ell^2)\otimes_{\max}C,
\]
then, for every Hilbert space $\cH$, without any
separability assumption,
\[
 B(\cH)\otimes_{\min}C
 =
 B(\cH)\otimes_{\max}C.
\]
\end{lemma}

\begin{proof}
The case $\cH=\{0\}$ is immediate, so assume that
$\cH\ne\{0\}$. Let
\[
 x=\sum_{j=1}^m T_j\otimes c_j
 \in B(\cH)\otimes C,
\]
where the tensor product in this display is algebraic, and put
\[
 \mathcal{B}_0=C^*(I_{\cH},T_1,\ldots,T_m)
 \subseteq B(\cH).
\]
Since $\mathcal{B}_0$ is separable, it admits a faithful
unital $*$-representation
\[
 \iota:\mathcal{B}_0\longrightarrow B(\ell^2).
\]
The inverse map
\[
 \iota^{-1}:\iota(\mathcal{B}_0)\longrightarrow
 \mathcal{B}_0\subseteq B(\cH)
\]
is a unital $*$-homomorphism. By Arveson's extension
theorem~\cite[Theorem~7.5]{Pa}, it extends to a unital
completely positive map
\[
 \Phi:B(\ell^2)\longrightarrow B(\cH).
\]
In particular, $\Phi(\iota(b))=b$ for every
$b\in\mathcal{B}_0$.

Set
\[
 y=\sum_{j=1}^m\iota(T_j)\otimes c_j
 \in B(\ell^2)\otimes C.
\]
Functoriality of the maximal tensor product for completely
positive contractions~\cite[Remark~12]{Pi96} implies that
\[
 \Phi\otimes\id_C:
 B(\ell^2)\otimes_{\max}C
 \longrightarrow B(\cH)\otimes_{\max}C
\]
is a unital completely positive map, and hence is contractive.
Since $(\Phi\otimes\id_C)(y)=x$, we obtain
\[
 \begin{aligned}
 \|x\|_{B(\cH)\otimes_{\max}C}
 &\le \|y\|_{B(\ell^2)\otimes_{\max}C}\\
 &=\|y\|_{B(\ell^2)\otimes_{\min}C}\\
 &=\|x\|_{B(\cH)\otimes_{\min}C}.
 \end{aligned}
\]
The first equality follows from the hypothesis.
The last equality follows from injectivity
of the minimal tensor product, applied to the faithful maps
$\mathcal{B}_0\hookrightarrow B(\cH)$ and
$\iota:\mathcal{B}_0\hookrightarrow B(\ell^2)$.
Since the reverse norm inequality always holds, the minimal
and maximal norms agree on the algebraic tensor product.
Passing to the completions proves the assertion.
\end{proof}

\begin{proof}[Proof of Corollary~\ref{t:main+amenableY}]
Put
\[
 C=C_u^*(\algA)
   =C^*(F_r)*\cA_1*\cdots *\cA_s,
 \qquad D=C\otm C^*(\groupS).
\]
The finite-dimensional factors $\cA_j$ are nuclear.
Thus Pisier's theorems for free groups and unital free
products~\cite[Theorems~0.1 and~0.2]{Pi96} give
\[
 B(\ell^2) \otimes_{\min}C
 =
 B(\ell^2)\otimes_{\max}C.\]
Applying Lemma~\ref{lem:min-max-arbitrary-H} gives
\begin{equation}\label{e:amenable-max-min}
 B(\cH)\otimes_{\min}C=B(\cH)\otm C.
\end{equation}

Since $\groupS$ is amenable, $C^*(\groupS)$ is nuclear.
In particular,
\[
 D=C\otimes_{\min}C^*(\groupS).
\]
Using associativity, nuclearity, and
\eqref{e:amenable-max-min}, we obtain
\[
\begin{aligned}
 B(\cH)\otm D
 &=(B(\cH)\otm C)\otm C^*(\groupS)\\
 &=(B(\cH)\otimes_{\min}C)
       \otimes_{\min}C^*(\groupS)\\
 &=B(\cH)\otimes_{\min}
       \bigl(C\otimes_{\min}C^*(\groupS)\bigr)\\
 &=B(\cH)\otimes_{\min}D.
\end{aligned}
\]

By \eqref{eq:positive-evaluations}, there is an
$\varepsilon>0$ such that
\[
 A\succeq\varepsilon I
 \quad\text{in}\quad B(\cH)\otimes_{\min}D.
\]
The preceding equality therefore gives
\[
 A\succeq\varepsilon I
 \quad\text{in}\quad
 B(\cH)\otm C_u^*(\algA)\otm C^*(\groupS).
\]
Theorem~\ref{t:main+maxY} now yields the asserted
factorization, including the half-extent bound and the
support restriction $u\in\lly_{\le M'}$ for some $M'\ge M$.
The target space can be taken to be a finite ampliation
of $\cH$, and is therefore separable whenever $\cH$ is
separable. If $\cH$ is infinite-dimensional, this
ampliation is unitarily isomorphic to $\cH$.
\end{proof}

\subsection{Limitations of spatial positivity}

The conclusion of Remark~\ref{r:matrix-coefficients} does not, in
general, extend to $B(\cH)$-coefficients under the same
uniform strict positivity hypothesis. 

\begin{example}
\label{r:operator-coefficients-obstruction}
\rm
Let $\cH=\ell^2$ and
$\Gamma=\operatorname{SL}_3(\Z)$, and choose a finite
symmetric generating set for $\Gamma$, so that its
generated semigroup is $\Gamma$ itself.

By \cite[Corollary~B]{ISW20}, $C^*(\Gamma)$ does not have
the local lifting property. Kirchberg's tensor-product
characterization of this property, stated in
\cite[Theorem~1.3]{Oz01}, therefore gives
\[
 X\in B(\cH)\otimes\C[\Gamma],
 \qquad \|X\|_{\min}<\|X\|_{\max}.
\]
Here the norms are induced by
$B(\cH)\otimes_{\min}C^*(\Gamma)$ and
$B(\cH)\otimes_{\max}C^*(\Gamma)$, respectively.
The element $X$ may be chosen with group-algebra
coefficients by approximating its finitely many
$C^*(\Gamma)$-coefficients while preserving the strict
norm inequality.

Choose
\[
 \|X\|_{\min}^2<t<\|X\|_{\max}^2,
 \qquad A_0=tI-X^*X,
\]
and let $A$ be obtained from $A_0$ by inserting the identity
in the $\algA$ factor. Its finite group support satisfies
the support condition of Corollary~\ref{t:main+amenableY}
for some
$M$, and its $\algA$-extent is zero. For every unital
$*$-representation $\pi$ of $\algA\otimes\C[\Gamma]$,
\[
 (\id_{B(\cH)}\otimes\pi)(A)
 \succeq\bigl(t-\|X\|_{\min}^2\bigr)I.
\]
However, $A_0$ is not positive in
$B(\cH)\otimes_{\max}C^*(\Gamma)$, since positivity there
would imply $\|X\|_{\max}^2\le t$.

Suppose that $A=B^*B$ for a finite polynomial
\[
 B=\sum_{q,u}B_{q,u}\otimes(q\otimes u),
 \qquad B_{q,u}\in B(\cH,\cE),
\]
where the $q$ are reduced words. Apply
$\id\otimes\tau_{\algA}\otimes\id$, where
$\tau_{\algA}(p^*q)=\delta_{p,q}$ by
Lemma~\ref{lem:orthonormal-reduced-words}. This gives
\[
 A_0=\sum_{u,v}Q_{u,v}\otimes u^{-1}v,
 \qquad Q_{u,v}=\sum_q B_{q,u}^*B_{q,v}.
\]
The finite matrix $Q=[Q_{u,v}]_{u,v}$ is positive over
$B(\cH)$. Therefore, with
$V=(I_\cH\otimes u)_u$ viewed as a column,
\[
 A_0=V^*(Q\otimes1)V\succeq0
 \quad\text{in }B(\cH)\otimes_{\max}C^*(\Gamma),
\]
a contradiction. Thus the obstruction persists even if
no extent bound is imposed on the factor.
\end{example}

We conclude with the following example that shows the uniform lower bound in
Corollary~\ref{t:main+amenableY} cannot be replaced by
strict positivity of the quadratic form in every spatial
evaluation.

\begin{example}\rm
\label{ex:below} 
Take $\algA=\C$, $\groupS=\Z^3$, and
$\cH=\ell^2(\N)$.
Let $p\geq0$ be a trigonometric polynomial on
$\mathbb T^3$ that is not a sum of hermitian squares.
Put
\[
 D=\operatorname{diag}(1,1/2,1/3,\ldots),
 \qquad A=I_{\cH}\otimes p+D\otimes1.
\]
For every unital $*$-representation
$\pi:\C[\Z^3]\to B(\cK)$ and every nonzero
$\xi\in\cH\otimes\cK$,
\[
 \bigl\langle
   (\id_{B(\cH)}\otimes\pi)(A)\xi,\xi
 \bigr\rangle
 \geq
 \bigl\langle(D\otimes I_{\cK})\xi,\xi\bigr\rangle
 >0.
\]
Nevertheless, $A$ admits no factorization $A=B^*B$
with $B$ an operator-valued polynomial.

Indeed, suppose that $A=B^*B$, where
\[
 B=\sum_{g\in F}B_g\otimes g,
 \qquad B_g\in B(\cH,\cE),
\]
and $F\subseteq\Z^3$ is finite.
Choose a state $\varphi$ on $B(\cH)$ annihilating compact operators. Then $\varphi(D)=0$, and applying $\varphi$ coefficientwise gives
\[
p=\sum_{g,h\in F}\varphi(B_g^*B_h)\,g^{-1}h.
\]
Since
\[
\bigl[\varphi(B_g^*B_h)\bigr]_{g,h\in F}\succeq0,
\]
this expresses $p$ as a finite sum of hermitian squares, a contradiction.
\qed
\end{example}

\appendix

\section{An elementary Følner-set proof of Corollary~\ref{t:main+amenableY}}

We also give a direct
F\o lner-set proof of Corollary~\ref{t:main+amenableY}.

\begin{proof}

Choose $\varepsilon>0$ such that
\[
 (\id_{B(\cH)}\otimes\pi)(A)\succeq\varepsilon I
\]
for every unital $*$-representation $\pi$ of
$\algA\otimes\C[\groupS]$; such an $\varepsilon$
exists by \eqref{eq:positive-evaluations}.

Choose a finite symmetric set $G\subseteq\groupS$ containing
$\idY$ and the group support of $A$, and write
\[
 A=\sum_{g\in G}A_g\otimes g,
 \qquad
 A_g=\sum_{a\in\cW_{2\vd}}A_{a,g}\otimes a.
\]
Coefficients outside the original support are understood
to be zero. Since $A=A^*$, we have $A_{g^{-1}}=A_g^*$.
For $a\in\algA$, let
\[
 \|a\|_u=\sup_\sigma\|\sigma(a)\|,
\]
where the supremum is over all unital $*$-representations
of $\algA$. This norm is finite, since the cyclic generators
are represented by unitaries and the representations of
the finite-dimensional factors are contractive. Put
\[
 C_0=\sum_{g\in G}\sum_{a\in\cW_{2\vd}}
       \|A_{a,g}\|\,\|a\|_u,
\]
and choose $0<\delta<1$ so that
\[
 \frac{\delta C_0}{1-\delta}<\frac{\varepsilon}{2}.
\]

Since $\groupS$ is amenable, there is a finite nonempty
right F\o lner set $F\subseteq\groupS$ such that
\[
 k_F(g):=
 \frac{\#\{(u,v)\in F^2:u^{-1}v=g\}}{|F|}
 =\frac{|F\cap Fg^{-1}|}{|F|}
 >1-\delta,
 \qquad g\in G.
\]
Notice that $k_F(g^{-1})=k_F(g)$ and $0\le k_F(g)\le1$.
The common-left-denominator argument in the proof of
Lemma~\ref{l:inf-to-rep-plus} gives $c\in\lly$ such that
$cF\subseteq\lly$. Since
\[
 (cu)^{-1}(cv)=u^{-1}v,
\]
replacing $F$ by $cF$ leaves every $k_F(g)$ unchanged.
We may therefore assume that
\[
 F\subseteq\lly_{\le M'}
\]
for some integer $M'\ge M$.

Define
\[
 A'_g=
 \begin{cases}
  k_F(g)^{-1}A_g,&g\in G,\\
  0,&g\notin G,
 \end{cases}
 \qquad
 A'=\sum_{g\in G}A'_g\otimes g.
\]
The symmetry of $G$ and the identity
$k_F(g^{-1})=k_F(g)$ imply that $A'=A'^*$.
For every unital $*$-representation $\pi$ of
$\algA\otimes\C[\groupS]$,
\begin{align*}
 \bigl\|(\id_{B(\cH)}\otimes\pi)(A'-A)\bigr\|
 &\le
 \sum_{g\in G}\sum_{a\in\cW_{2\vd}}
 \left|\frac1{k_F(g)}-1\right|
 \|A_{a,g}\|\,\|a\|_u\\
 &\le\frac{\delta C_0}{1-\delta}
 <\frac{\varepsilon}{2}.
\end{align*}
Consequently,
\begin{equation}
\label{e:amenable-corrected-positive}
 (\id_{B(\cH)}\otimes\pi)(A')
 \succeq\frac{\varepsilon}{2}I.
\end{equation}

Consider the finite operator matrix
\[
 T=\bigl[A'_{u^{-1}v}\bigr]_{u,v\in F}
 \in B\bigl(\cH\otimes\ell^2(F)\bigr)
       \otimes\algA_{2\vd}.
\]
We claim that $T$ is positive on every finite-dimensional
unital $*$-representation of $\algA$. To see this, let
$\sigma:\algA\to B(\cK)$ be such a representation and let
$R$ be the right regular representation of $\groupS$,
\[
 R_g\delta_v=\delta_{vg^{-1}}.
\]
Its $(u,v)$ matrix entry is $1$ precisely when $g=u^{-1}v$.
Apply \eqref{e:amenable-corrected-positive} to the
representation
\[
 a\otimes g\longmapsto\sigma(a)\otimes R_g,
\]
and compress the $\ell^2(\groupS)$ factor to $\ell^2(F)$.
After the canonical permutation of tensor factors, the
result is
\[
 \bigl(\id_{B(\cH\otimes\ell^2(F))}\otimes\sigma\bigr)(T)
 \succeq\frac{\varepsilon}{2}I.
\]
This proves the claim.

By Theorem~\ref{thm:main}, there exist a Hilbert space
$\cE$ and an element
\[
 D\in
 B\bigl(\cH\otimes\ell^2(F),\cE\bigr)
 \otimes\algA_\vd
 \qquad\text{such that}\qquad D^*D=T.
\]
We may take
\[
 \cE=
 \bigl(\cH\otimes\ell^2(F)\bigr)
 \otimes\C^{N(\vd)}.
\]
For $u\in F$, let
$J_u:\cH\to\cH\otimes\ell^2(F)$ be the isometry
$J_u h=h\otimes\delta_u$, and set
\[
 D_u=D(J_u\otimes\idA)
 \in B(\cH,\cE)\otimes\algA_\vd.
\]
Then $D_u^*D_v=A'_{u^{-1}v}$. Define
\[
 B=\frac1{\sqrt{|F|}}\sum_{u\in F}D_u\otimes u.
\]
Expanding each $D_u$ in the basis $\cW_\vd$ shows that $B$
has the asserted support, with zero coefficients for
$u\in\lly_{\le M'}\setminus F$. Moreover,
\begin{align*}
 B^*B
 &=\frac1{|F|}\sum_{u,v\in F}
       D_u^*D_v\otimes u^{-1}v\\
 &=\frac1{|F|}\sum_{u,v\in F}
       A'_{u^{-1}v}\otimes u^{-1}v\\
 &=\sum_{g\in G}k_F(g)A'_g\otimes g
 =A.
\end{align*}
The displayed choice of $\cE$ is a finite ampliation of
$\cH$, and is therefore separable whenever $\cH$ is
separable. If $\cH$ is infinite-dimensional,
choose a unitary $U:\cE\to\cH$ and replace $B$ by
$(U\otimes\idA\otimes\idY)B$ to obtain coefficients in
$B(\cH)$ without changing $B^*B$.
\end{proof}

\end{document}